\documentclass[11pt, reqno]{amsart}
\usepackage{enumitem}
\usepackage{amsthm,amsmath,amssymb}

\usepackage{mathrsfs}

\numberwithin{equation}{section}

\newtheorem{theorem}{Theorem}[section]
\newtheorem{definition}[theorem]{Definition}
\theoremstyle{plain}
\newtheorem{proposition}[theorem]{Proposition}
\newtheorem{lemma}[theorem]{Lemma}
\newtheorem{corollary}[theorem]{Corollary}
\newtheorem{remark}[theorem]{Remark}

\def\be{\begin{equation}}
	\def\ee{\end{equation}}

\def\d{\nabla}

\DeclareMathOperator{\Ric}{Ric}

\DeclareMathOperator{\supp}{supp}
\DeclareMathOperator{\vol}{Vol}

\DeclareMathOperator{\Hess}{Hess}

\numberwithin{equation}{section}

\begin{document}
\title[]{Heat kernel estimate for the Laplace-Beltrami operator under Bakry-\'Emery Ricci curvature condition and applications}
\author{ Xingyu Song, Ling Wu, and Meng Zhu}
\address{School of Mathematical Sciences and Shanghai Key Laboratory of PMMP, East China Normal University, Shanghai 200241, China}
\email{52215500013@stu.ecnu.edu.cn, 52215500012@stu.ecnu.edu.cn, mzhu@math.ecnu.edu.cn}
\date{}

\begin{abstract}
	We establish a Gaussian upper bound of the heat kernel for the Laplace-Beltrami operator on complete Riemannian manifolds with Bakry-\'Emery Ricci curvature bounded below. As applications, we first prove an $L^1$-Liouville property for non-negative subharmonic functions when the potential function of the  Bakry-\'Emery Ricci curvature tensor is of at most quadratic growth. Then we derive lower bounds of the eigenvalues of the Laplace-Beltrami operator on closed manifolds. An upper bound of the bottom spectrum is also obtained.

\end{abstract}
\maketitle

\section{Introduction}
Let $(M^{n},g)$ be an $n$-dimensional Riemmanian manifold. The Bakry-\'Emery Ricci curvature tensor of $M$ (\cite{BE}) is defined as
\be\label{BE curvature}
\Ric_f:=\Ric+\Hess f,
\ee
where $f$ is a smooth function on $M$ (called the potential function), and $\Ric$ and $\Hess f$ denote the Ricci curvature tensor  and the hessian of $f$, respectively. It is clear that when $f$ is a constant, $\Ric_f$ reduces to the Ricci curvature tensor. Also, manifolds with lower Bakry-\'Emery Ricci curvature bound are closely related to the singularity analysis of the Ricci flow, Ricci limit spaces, and stationary black holes (see e.g., \cite{Ham, Per, St1, St2, LV, GKW}). Therefore, many efforts have been made in extending the results under Ricci curvature condition to Bakry-\'Emery Ricci curvature condition.

Since $\Ric_f$ appears in a Bochner type formula for the weighted Laplace operator $\Delta_f:= \Delta-\left<\d f,\d\right>$ (see \cite{BE}) which is self-adjoint with respect to the weighted measure $e^{-f} dv$ on $M$, the Bakry-\'Emery Ricci curvature can be considered as the "Ricci curvature" for the smooth metric measure space $(M^n, g, e^{-f}dv)$. Here $\Delta$ is the Laplace-Beltrami operator, $\d f$ is the gradient of $f$, and $dv$ is the volume form of $g$. Thus, while doing analysis under Bakry-\'Emery Ricci curvature condition, one would
naturally consider replacing $\Delta$ and $dv$ with the weighted Laplace operator $\Delta_f$ and volume form $e^{-f}dv$ (see e.g., \cite{BE, Lo, WeWy, FLZ, MW, Liu, WZ, Wu2, WuWu, WuWu2}).


Let us mention that there are many similar studies under assumptions on the $m$-Bakry-\'Emery Ricci curvature, $\Ric^m_f=\Ric+\Hess f-\frac{\d f\otimes \d f}{m-n}$, where $m\in(n,+\infty)$ is a constant (see e.g., \cite{BQ,  Lix1, Liy, SuZh, Wu1}). The Bakry-\'Emery Ricci curvature corresponds to the case $m=+\infty$, and hence it is also called the $\infty$-Bakry-\'Emery Ricci curvature.

In this paper, we treat the Bakry-\'Emery Ricci curvature as a generalization of the Ricci curvature on $(M^n, g)$ with the original volume measure $dv$, and study the properties of elliptic and the heat equations related to the Laplace-Beltrami operator. This type of study has previously been conducted by Q. Zhang and the third author, and was applied in the extension of the Cheeger-Colding-Naber theory and the proof of compactness theorems for gradient Ricci solitons (see \cite{ZZ1}, \cite{ZZ2}, \cite{GPSS}, \cite{LLW}).

First,  we investigate estimates of the heat kernel. On manifolds with Ricci curvature bounded below, there have already been a long history and many classical results about the heat kernel estimates, the readers may refer to \cite{Da, Gr, SC} and the references therein. For manifolds with Bakry-\'Emery Ricci curvature, or $m$-Bakry-\'Emery Ricci curvature bounded below, the heat kernel estimates for the weighted Laplace operator $\Delta_f$ have been studied in many literatures recently (see \cite{Qian, Lix1, Lix2, Wu3, WuWu, WuWu2}).

Our first main result is a Gaussian upper bound of the heat kernel of $\Delta$ depending on the lower bound of the Bakry-\'Emery Ricci curvature and the bound of  the potential function $f$.

\begin{theorem}(Gaussian upper bound of the heat kernel)\label{Gaussian}
	Let $(M^n,g)$ be a complete Riemannian manifold and $H(x,y,t)$ the heat kernel of $\Delta$ on $M\times M\times (0,+\infty)$. Pick a fixed point $o\in M$ and $R>0$. Suppose that $\Ric_f\ge-Kg$ on $B_{3R}(o)$ for some constant $K\ge0$. Then for any $\epsilon>0$, there exist constants $C_1(n,\epsilon)$ and $C_2(n)$, such that
	\begin{equation}\label{4.8}
		H(x,y,t)\le\frac{C_1(n,\epsilon)e^{C_2(n)(Kt+L(R))}}{\vol(B_{\sqrt{t}}(x))^{\frac{1}{2}}\vol(B_{\sqrt{t}}(y))^{\frac{1}{2}}}e^{\left(-\frac{d^2(x,y)}{(4+\epsilon)t}\right)}
	\end{equation}	
	for all $x,y \in B_{\frac{1}{2}R}(o)$ and $0<t<R^2/4$, where $L(R)=\sup\limits_{B_{3R}(o)}|f|$ and $\lim\limits_{\epsilon \to 0} C_1(n,\epsilon)=\infty$.
\end{theorem}	

\begin{remark}
Note that in general the heat kernel $H(x,y,t)$ of the Laplace-Beltrami operator $\Delta$ and the heat kernel $H_f(x,y,t)$ of the weighted Laplace operator $\Delta_f$ are not equivalent. Hence, one cannot expect to derive \eqref{4.8} directly from the estimates of  $H_f(x,y,t)$. Indeed, on the Euclidean space $(\mathbb{R}^n, g_0)$ with $f=x_1$, we have (see e.g. \cite{WuWu})  
$$H_f(x,y,t)=\frac{e^{\frac{2(x_1+y_1)-t}{4}}}{(4\pi t)^{-\frac{n}{2}}}e^{-\frac{|x-y|^2}{4t}}=e^{\frac{2(x_1+y_1)-t}{4}}H(x,y,t).$$

\end{remark}
	


The idea of the proof is standard. A key ingredient is a relative volume comparison theorem (Theorem \ref{volume element comparison} below) under the assumptions in the theorem above. Using the volume comparison, we can derive a Sobolev inequality, from which a parabolic mean value property follows. Then the mean value property and Davies' double integral estimate will imply the heat kernel upper bound.

With the help of the heat kernel upper bound, one can prove an $L^1$-Liouville property for non-negative subharmonic functions. On manifolds with Ricci curvature bounded below, this was first proved by Li-Schoen \cite{LiSc} and Li \cite{Li} , where the lower bound of the Ricci curvature is allowed to have certain decay. We show that the same result still holds on manifolds with Bakry-\'Emery Ricci curvature bounded below.


\begin{theorem}\label{Liouville theorem}
	Let $(M^n,g)$ be a complete  noncompact Riemannian manifold with $\Ric_f\ge -Kg$ for some constant $K\ge0$. Assume that there exist non-negative constants $a$ and $b$ such that
	\begin{equation}\label{gradient linear growth}	
		|f(x)|\le ar^2(x)+b \ for \ all \ x\in M,
	\end{equation}	
	where $r(x)=d(x,o)$ is the geodesic distance function to a fixed point $o\in M$. Then any non-negative $L^1$-integrable subharmonic function on $M$ must be identically constant. In particular, any $L^1$-integrable harmonic function must be identically constant.	
\end{theorem}

\begin{remark}
	In \cite{Yau1}, Yau obtained the $L^{\infty}$ Liouville property of harmonic functions on manifolds with nonnegative Ricci curvature. Later, Yau \cite{Yau2} showed that for $1<p<+\infty$ the $L^{p}$ Liouville property for subharmonic functions actually holds without any curvature assumption.
	
	The $L^1$ Liouville property for subharmonic functions with respect to the weighted Laplace operator $\Delta_f$ on manifolds with Bakry-\'Emery Ricci curvature (or $m$-Bakry-\'Emery Ricci curvature, resp.) bounded below was proved by Wu-Wu \cite{WuWu2} (X. Li \cite{Lix1}, resp.).
\end{remark}

The conditions in the above theorem are especially satisfied on so called gradient Ricci solitons. A gradient Ricci soliton is a Riemannian manifold $(M^n, g)$ with constant Bakry-\'Emery Ricci curvature, namely,
\be\label{gradient Ricci soliton}
\Ric + \Hess f=\lambda g
\ee
for some constant $\lambda$. It is called a shrinking, steady, or expanding Ricci soliton when $\lambda>0$, $=0$, or $<0$, respectively.

For gradient Ricci solitons, it is well known that
\be\label{gradient f}
S+|\d f|^2=2\lambda f + C,
\ee
where $S$ is the scalar curvature of $M$. Then from  the results in  \cite{Chen} and \cite{CZ} , we know that $|f|$ is at most of quadratic and linear growth on gradient shrinking and steady Ricci solitons, respectively. This implies that \eqref{gradient linear growth} holds.


For gradient expanding solitons, it is showed in \cite{PRS} that $S\geq n\lambda$, then \eqref{gradient f} implies that $\left|\d \sqrt{-f-\frac{C-n\lambda}{2\lambda}}\right|\leq \sqrt{-\frac{\lambda}{2}}$. Hence
$$\sqrt{-f(x)-\frac{C-n\lambda}{2\lambda}}\leq \sqrt{-\frac{\lambda}{2}}r(x)+\sqrt{-f(o)-\frac{C-n\lambda}{2\lambda}},$$
and the conditions in Theorem \ref{Liouville theorem} continue to hold with $f$ replaced by $\tilde{f}=f(x)+\frac{C-n\lambda}{2\lambda}$.

Therefore, Theorem \ref{Liouville theorem} implies the following Liouville theorem for complete gradient Ricci solitons.

\begin{theorem}\label{Liouvile theorem for Ricci soliton}
	Let $(M^n,g)$ be a complete gradient Ricci soliton. Then any non-negative $L^1$-integrable subharmonic function on $M$ must be identically constant. In particular, any $L^1$-integrable harmonic function on $M$ must be identically constant.
\end{theorem}


\begin{remark}
	In \cite{MuSe}, Munteanu-Sesum proved that on a gradient shrinking K\"ahler-Ricci soliton or a gradient steady Ricci soliton, harmonic functions with finite energy must be constant.
\end{remark}

By Theorem \ref{Liouville theorem}, following the arguments of Li in \cite{Li} we can prove a uniqueness theorem for $L^1$-solutions of the heat equation.

\begin{theorem}\label{L^1}
	Let $(M^n,g)$ be a complete  noncompact Riemannian manifold with $\Ric_f\ge -Kg$ for some constant $K\ge0$. Assume that  there exist non-negative constants $a$ and $b$ such that
	\begin{equation}	
		|f|(x)\le ar^2(x)+b \ for \ all \ x\in M.\nonumber
	\end{equation}	
	If $u(x,t)$ is a non-negative function defined on $M\times [0,+\infty)$ satisfying
	\begin{equation}
		\left(\Delta-\partial_t\right)u(x,t)\ge 0,\ \int_M u(x,t)dv<+\infty \nonumber
	\end{equation}
	for all $t>0$, and $\lim\limits_{t\rightarrow 0}\int_M u(x,t)=0$, then $u(x,t)\equiv 0$.
	
	In particular, any $L^1$-solution of the heat equation is uniquely determined by its initial data in $L^1$.
\end{theorem}
Another application of the heat kernel estimate is that one can get Li-Yau type \cite{LY} lower bound estimates for the eigenvalues of $\Delta$ on closed manifolds. More precisely, we show that
\begin{theorem}\label{Eigenvalue Estimates theorem}
	Let $(M^n,g)$ be a closed Riemannian manifold with $\Ric_f\ge-Kg$ and $|f|\le L$. Let $0=\lambda_0<\lambda_1\leq \lambda_2\leq\cdots\leq \lambda_k\leq \cdots$ be the eigenvalues of the Laplace-Beltrami operator. Then there exist constants $C_3(n)$ and $C_4(n)$ such that
	\begin{equation}\label{eigenvalue estimates}
		\lambda_k\ge  \frac{C_3(n)(k+1)^{\frac{2}{n}}}{D^2}e^{-C_4(n)(KD^2+L)}
	\end{equation}	
	for all $k\ge1$, where $D$ is an upper bound of the diameter of $M$.
\end{theorem}

Finally, to get an upper bound of $\lambda_1$, we may relax the requirements on the manifold and $f$. Indeed, we can obtain a Cheng type upper bound of the bottom spectrum $\mu_1(\Delta)$ of $\Delta$ on complete manifolds with $f$ being of at most linear growth.

\begin{theorem}\label{upper bound of spectrum}
		Let $(M^n,g)$ be a complete Riemannian manifold with $\Ric_f\ge -(n-1)Kg$ for some constant $K\ge0$. Assume that there exist non-negative constants $\tilde{a}$ and $\tilde{b}$ such that 
	\begin{equation}
		|f|(x)\le \tilde{a}r(x,o)+\tilde{b}\ for\ all\ x\in M.\nonumber
	\end{equation}
Then we have 
\begin{equation}
	\mu_1(\Delta)\le\frac{1}{4}\left(2\tilde{a}+(n-1)\sqrt{K}\right)^2.\nonumber
\end{equation}
In particular, if $f$ is of sublinear growth, then the bottom spectrum of the Laplacian has the following sharp upper bound:
\begin{equation}
	\mu_1(\Delta)\le\frac{1}{4}\left(n-1\right)^2K.\nonumber
\end{equation}
\end{theorem}
When $f$ is constant, i.e, $\Ric\ge-(n-1)Kg$, Theorem \ref{upper bound of spectrum} reduces to Cheng's theorem (see Theorem 4.2 in \cite{CSY}). This incidentally indicates that our estimate is sharp.\\

The rest of the paper is organized as follows. In section 2, we derive a Laplacian and volume comparison theorem for complete Riemannian manifolds with Bakry-\'Emery Ricci curvature bounded below. Using the comparison theorem, we get Poincar\'e and Sobolev inequalities.  Section 3 is devoted to the proofs of a mean value inequality for non-negative subsolutions of the heat equation and the Gaussian upper bound \eqref{4.8} of the heat kernel. The $L^1$-Liouville property for non-negative subharmonic functions and the uniqueness theorem for $L^1$-solutions of the heat equation will be obtained in section 4. Then for the purpose of comparison, we take a detour  in section 5 to discuss an $L^{\infty}$-Liouville property for harmonic function with polynomial growth.  Finally, lower bounds of the eigenvalues and an upper bound of the bottom spectrum of the Beltrami-Laplace operator are shown in section 6.

\section{Poincar\'e and Sobolev inequalities  }

In this section, we derive local Poincar\'e and Sobolev inequalities which are instrumental in proving the main results of the paper.
A crucial step is to show the following volume comparison theorem, the idea of whose proof is similar to \cite{Yang}, \cite{MW} and \cite{ZZ1}, but our result is slightly more general.

For a fixed point $o\in M$ and $R>0$, we define
\begin{equation}\label{eq2}
	L(R)=\sup\limits_{x\in B_{3R}(o)}|f(x)|,
\end{equation}
where $ B_{3R}(o)$ is the geodesic ball centered at $o\in M$ with radius $3R$.

\begin{theorem}\label{volume element comparison}
Let $(M^n, g)$ be a complete Riemannian manifold with $\Ric_f\ge-Kg$ for some constant $K\ge0$. Then the following conclusions are true.
\\
(a)(Laplacian comparison) Let $r=d(y,p)$ be the distance from any point $y$ to some fixed point $p\in B_{R}(o)$ with $0<r<R$. Then for $0<r_1<r_2<R$, we have
\begin{equation}\label{Laplacian}
	\int_{r_1}^{r_2}(\Delta r-\frac{n-1}{r})dr \le \frac{K}{6}(r_2^2-r_1^2)+6L(R).
\end{equation} 
\\
(b)(Volume element comparison)Take any point $p\in B_{R}(o)$ and denote the volume form in geodesic polar coordinates centered at $p$ with $J(r,\theta,p)drd\theta$, where $r>0$ and $\theta\in S_p(M)$, a unit tangent vector at $p$. Then for $0<r_1<r_2<R$, we have
\begin{equation}\label{AC}
	\frac{J(r_2,\theta,p)}{J(r_1,\theta,p)}\le \left(\frac{r_2}{r_1}\right)^{n-1} e^{\frac{K}{6}(r_2^2-r_1^2)+6L(R)} .
\end{equation}
\\
(c)(Volume comparison)For any $p\in B_{R}(o),\ 0<r_1<r_2<R$, we have
\begin{equation} \label{VC}
	\frac{\vol(B_{r_2}(p))}{\vol(B_{r_1}(p))}\le \left(\frac{r_2}{r_1}\right)^{n} e^{\frac{K}{6}(r_2^2-r_1^2)+6L(R)},
\end{equation}
where $\vol(.)$ denotes the volume of a region.
\end{theorem}
\textit{Proof of part (a)} Let $r=d(y,p)$ be the distance from any point $y$ to some fixed point $p\in B_{R}(o)$ with $0<r<R$ and $\gamma:[0,r]\to M$ a normal minimal geodesic with $\gamma(0)=p$ and $\gamma(r)=y$. Then we know $\gamma(t)  \subset B_{3R}(o)$ and $y\in B_{3R}(o)$.\\
From the Bochner formula, we have 
\begin{equation}
0=\frac{1}{2}\Delta|\nabla r|^2=|\Hess r|^2+\left<\nabla\Delta r,\nabla r\right>+\Ric(\partial r,\partial r). \nonumber
\end{equation}
By using Cauchy-Schwarz inequality $|\Hess r|^2\ge\frac{(\Delta r)^2}{n-1}$, it yields
\begin{equation}
\frac{\partial}{\partial r}(\Delta r)+\frac{(\Delta r)^2}{n-1}\le -\Ric(\partial r,\partial r),\nonumber
\end{equation}
which is equivalent to 
\begin{equation}\label{eqq2.5}
\frac{1}{r^2}\frac{\partial}{\partial r}(r^2\Delta r)+\frac{1}{n-1}\left(\Delta r-\frac{n-1}{r}\right)^2\le\frac{n-1}{r^2}-\Ric(\partial r,\partial r).
\end{equation}
Multiplying both sides of \eqref{eqq2.5} by $r^2$ and integrating from $0$ to $r$, we get 
\begin{equation}\label{eqq2.6}
\Delta r \le \frac{n-1}{r}-\frac{1}{r^2}\int_{0}^{r}t^2\Ric(\gamma'(t),\gamma'(t))dt.
\end{equation}
We observe that 
\begin{equation}
\begin{aligned}
	\Ric(\gamma'(t),\gamma'(t)) &\ge-K-\Hess f(\gamma'(t),\gamma'(t))\\
	&=-K-\left<\nabla_{\gamma'(t)}\nabla f,\gamma'(t)\right>\\
	&=-K-\gamma'(t)\left<\nabla f,\gamma'(t)\right>+\left<\nabla f,\nabla_{\gamma'(t)}\gamma'(t)\right>\\
	&=-K-f''(t), \nonumber
\end{aligned}
\end{equation}
where $f(t):=f(\gamma(t))$.\\
Hence, \eqref{eqq2.6} becomes 
\begin{equation}\label{2.7}
\begin{aligned}
	\Delta r-\frac{n-1}{r} &\le \frac{1}{r^2}\int_{0}^{r}\left[Kt^2+t^2f''(t)\right]dt\\
	&= \frac{1}{r^2}\left[\frac{K}{3}r^3+r^2f'(r)-2\int_{0}^{r}f'(t)tdt\right]\\
	&=\frac{K}{3}r+f'(r)-\frac{2}{r}f(r)+\frac{2}{r^2}\int_{0}^{r}f(t)dt.
\end{aligned}
\end{equation}
For $0<r_1<r_2<R$, integrating \eqref{2.7} from $r_1$ to $r_2$ yields
\begin{equation}
\begin{aligned}
	\int_{r_1}^{r_2}\left(\Delta r-\frac{n-1}{r}\right)dr &\le\frac{K}{6}(r_2^2-r_1^2)+\int_{r_1}^{r_2}\left(f'(r)-\frac{2}{r}f(r)\right)dr-2\int_{r_1}^{r_2}\int_{0}^{r}f(t)dtd\frac{1}{r}\\
	&=\frac{K}{6}(r_2^2-r_1^2)+f(r_2)-f(r_1)-2\int_{r_1}^{r_2}\frac{f(r)}{r}dr\\
	&\ \ \ -2\left[\left.\frac{1}{r}\int_{0}^{r}f(t)dt\right|^{r_2}_{r_1}-\int_{r_1}^{r_2}\frac{f(r)}{r}dr\right]\\
	&=\frac{K}{6}(r_2^2-r_1^2)+f(r_2)-f(r_1)-\frac{2\int_{0}^{r_2}f(t)dt}{r_2}+\frac{2\int_{0}^{r_1}f(t)dt}{r_1}\\
	&\le \frac{K}{6}(r_2^2-r_1^2)+6L(R).
\end{aligned}
\end{equation}

\textit{Proof of part (b)} By the first variation of the area 
\begin{equation}
\Delta r=\frac{J'(r,\theta,p)}{J(r,\theta,p)},\nonumber
\end{equation}
where $r(y)=d(y,p)$.
For $0<r_1<r_2<R$, integrating this from $r_1$ to $r_2$ we get 
\begin{equation}
\int_{r_1}^{r_2}\frac{\partial}{\partial r} \ln J(r,\theta,p)dr \le \int_{r_1}^{r_2}\frac{n-1}{r}dr+\frac{K}{6}(r_2^2-r_1^2)+6L(R).\nonumber
\end{equation}
Then we have 
\begin{equation}
\frac{J(r_2,\theta,p)}{J(r_1,\theta,p)} \le \left(\frac{r_2}{r_1}\right)^{n-1}e^{\frac{K}{6}(r_2^2-r_1^2)+6L(R)}.\nonumber
\end{equation}

\textit{Proof of part (c)} For $0<r_1<r_2<R$, we have
\begin{equation}
\frac{\vol (B_{r_2}(p))}{\vol( B_{r_1}(p))}=\frac{\int_{0}^{r_2}\int_{S^{n-1}}J(r,\theta,p)d\theta dr}{\int_{0}^{r_1}\int_{S^{n-1}}J(r,\theta,p)d\theta dr}=\frac{\frac{r_2}{r_1}\int_{0}^{r_1}\int_{S^{n-1}}J(\frac{r_2}{r_1}t,\theta,p)d\theta dt}{\int_{0}^{r_1}\int_{S^{n-1}}J(t,\theta,p)d\theta dt},\nonumber
\end{equation}
where $S^{n-1}$ denotes the unit sphere in $\mathbb{R}^n$ and $d\theta$ is its volume element.\\
By volume element comparison, for $0<t<r_1$, we have
\begin{equation}
	\begin{aligned}
		J\left(\frac{r_2}{r_1}t,\theta,p\right)&\le J(t,\theta,p)\left(\frac{r_2}{r_1}\right)^{n-1}e^{\frac{K}{6}\left(\frac{r_2^2}{r_1^2}t^2-t^2\right)+6L(R)}\\
		&\le J(t,\theta,p)\left(\frac{r_2}{r_1}\right)^{n-1}e^{\frac{K}{6}(r_2^2-r_1^2)+6L(R)}. \nonumber
	\end{aligned}
\end{equation}
Integrating this gives that
\begin{equation}
\frac{\vol (B_{r_2}(p))}{\vol( B_{r_1}(p))}\le\left(\frac{r_2}{r_1}\right)^{n}e^{\frac{K}{6}(r_2^2-r_1^2)+6L(R)} .\nonumber
\end{equation} 
\qed \\

By Theorem \ref{volume element comparison}, following  Buser's proof \cite{Bu} or Saloff-Coste's alternative proof (Theorem 5.6.5 in \cite{SC}), we can get a local Neumann Poincar\'e inequality, see also Lemma 3.1 in \cite{MW}.
\begin{lemma} \label{Poincare }
	Let $(M^n,g)$ be a complete Riemannian manifold with $\Ric_f\ge-Kg$ for some constant $K\ge0$. Then for any $p\in B_{R}(o)$, there exist constants $c_1$ and $c_2$ depending only on $n$ such that
	\begin{equation}
		\int_{B_{r}(p)} |u-u_{B_{r}(p)}|^2 \le c_1e^{c_2(Kr^2+L(R))}r^2\int_{B_{r}(p)} |\nabla u|^2 \nonumber
	\end{equation}
	for all $0<r<R$, where $u\in C^{\infty}\left(B_{r}(p)\right)$ and $u_{B_{r}(p)}=\frac{\int_{B_{r}(p)} u}{\vol (B_{r}(p))}.$
\end{lemma}

Combining Theorem \ref{volume element comparison} and Lemma \ref{Poincare }, and using a similar argument as in the proof of Lemma 3.2 in \cite{MW}, we obtain a local Neumann Sobolev inequality.

\begin{theorem} \label{Sobolev}
	Let $(M^n,g)$ be a complete Riemannian manifold with $\Ric_f\ge-Kg$ for some constant $K\ge0$. Then there exist constants $\mu=4n-2>2$, $c_3$ and $c_4$ depending only on $n$ such that
	\begin{equation}\label{Neumann Sobolev}
		\left(\int_{B_{r}(o)}|u-u_{B_{r}(o)}|^{\frac{2\mu}{\mu-2}} \right)^{\frac{\mu-2}{\mu}} \le \frac{c_3e^{c_4(Kr^2+L(R))}}{\vol (B_{r}(o))^\frac{2}{\mu}}r^2\int_{B_{r}(o)} |\nabla u|^2
	\end{equation}
	for all $0<r<R$, where $u\in C^{\infty}\left(B_{r}(o)\right)$ and $u_{B_{r}(o)}=\frac{\int_{B_{r}(o)} u}{\vol (B_{r}(o))}.$
\end{theorem}

By the Minkowski inequality and applying \eqref{Neumann Sobolev}, it is well known that one can get the following Sobolev inequality.
\begin{theorem} \label{Sobolev2}
	Let $(M^n,g)$ be a complete Riemannian manifold with $\Ric_f\ge-Kg$ for some constant $K\ge0$. Then there exist constants $\mu=4n-2>2$, $c_5$ and $c_6$, all depending only on $n$ such that
	\begin{equation}\label{eq10}
		\left(\int_{B_{r}(o)}|u|^{\frac{2\mu}{\mu-2}} \right)^{\frac{\mu-2}{\mu}} \le \frac{c_5e^{c_6(Kr^2+L(R))}}{\vol (B_{r}(o))^\frac{2}{\mu}}r^2\int_{B_{r}(o)} (|\nabla u|^2+r^{-2}u^2)
	\end{equation}
	for all $0<r<R$, where $u\in C^{\infty}\left(B_{r}(o)\right)$.
\end{theorem}

\section{  Mean value inequality and Gaussian upper bounds of the heat kernel }

In this section, we apply Theorem \ref{Sobolev2} to prove a mean value inequality for the non-negative subsolution of the heat equation  and the Gaussian upper bound \eqref{4.8} of the heat kernel on complete Riemannian manifolds with Bakry-\'Emery Ricci curvature bounded below.

In the following context, for function $u$ on $M$, the $L^q$ norm on a domain $\Omega\subset M$ is denoted by
\begin{equation}
	||u||_{q,\Omega}={\left(	\int_{\Omega}|f|^q\right)^{\frac{1}{q}}}.\nonumber
\end{equation}
$||u||_q$ denotes the $L^q$ norm of $u$ on $M$.\\

\begin{proposition}(Mean value inequality)\label{mean value inequality}
	Let $(M^n,g)$ be a complete Riemannian manifold with $\Ric_f\ge-Kg$ for some constant $K\ge0$. For any real number $s$ and any $0<\delta<\delta'\le1$, let $u$ be a smooth non-negative  subsolution of the heat equation in the cylinder $Q=B_{r}(o)\times(s-r^2,s)$, $0<r<R$.
	
For $2\le p <\infty$, there exist constants $\tilde{c}_1(n)$ and $\tilde{c}_2(n)$ such that
	\begin{equation}\label{4.1}
		\sup\limits_{Q_\delta}u^p\le \frac{\tilde{c}_1(n)e^{\tilde{c}_2(n)(Kr^2+L(R))}}{(\delta'-\delta)^{4n}r^2\vol (B_{r}(o))} \cdot\int_{Q_{\delta'}} u^pdvdt.
	\end{equation}

For $0<p<2$, there exist constants $\tilde{c}_3(n,p)$ and $\tilde{c}_4(n)$ such that
\begin{equation}\label{eqqu4.2}
	\sup\limits_{Q_\delta}u^p\le \frac{\tilde{c}_3(n,p)e^{\tilde{c}_4(n)(Kr^2+L(R))}}{(\delta'-\delta)^{4n}r^2\vol (B_{r}(o))} \cdot\int_{Q_{\delta'}} u^pdvdt.
\end{equation}
Here $Q_\delta=B_{\delta r}(o)\times(s-\delta r^2,s)$, $Q_{\delta'}=B_{\delta' r}(o)\times(s-\delta' r^2,s)$, $L(R)=\sup\limits_{B_{3R}(o)}|f|$.

	\proof The proof is similar to Theorem 5.2.9 in \cite{SC}.  We need to carefully examine the explicit coefficients of the mean value inequality in terms of the Sobolev constants in \eqref{eq10}.

	Without loss of generality we may assume $\delta'=1.$ Denote $B=B_{r}(o)$ for simplicity. For any non-negative function $ \phi \in C_0^\infty(B),$ we have

	\begin{equation}
		\int_{B}(\phi u_t+\left<\nabla u,\nabla \phi \right>)dv\le0.\nonumber
	\end{equation}
In particular, when $ \phi =\Phi^2u,\Phi \in C_0^\infty(B)$, we obtain
	
	\begin{equation}
		\begin{aligned}
			\int_{B}(\Phi^2 uu_t+\Phi^2|\nabla u|^2)dv&\le 2\left|\int_{B}u\Phi\left<\nabla u,\nabla \Phi \right> dv\right|\\
			&\le3\int_{B}|\nabla \Phi|^2u^2dv+\frac{1}{3}\int_{B}\Phi^2|\nabla u|^2dv.\nonumber
		\end{aligned}
	\end{equation}
It then implies that
	
	\begin{equation}
		\begin{aligned}
			\int_{B}(2\Phi^2 uu_t+|\nabla (\Phi u)|^2)dv&\le
			2\int_{B}\Phi^2  uu_tdv+\frac{4}{3}\int_{B}\Phi^2|\nabla u|^2dv+
			4\int_{B}|\nabla \Phi |^2  u^2dv\\
			&\le 10||\nabla \Phi ||_\infty^2\int_{\supp (\Phi)}u^2dv.\nonumber
		\end{aligned}
	\end{equation}
For any smooth non-negative function $\lambda(t)$ of the time variable $t$, which will be chosen later,  we get
	\begin{equation}
		\begin{aligned}
			\frac{\partial}{\partial t}	
			&\left(\int_{B}(\lambda \Phi u)^2dv\right)+\lambda^2\int_{B}|\nabla(\Phi  u)|^2dv\\
			&\le2\lambda|\lambda'|\sup\Phi^2\int_{\supp(\Phi)}u^2dv+\lambda^2\left(2\int_{B}\Phi^2  uu_tdv+\int_{B}|\nabla(\Phi  u)|^2dv\right)\\
			&\le C\lambda(\lambda||\nabla \Phi ||_\infty^2+|\lambda'|\sup\Phi^2)\int_{\supp(\Phi)}u^2dv,\nonumber
		\end{aligned}
	\end{equation}
where $C$ is a constant which will change from line to line in the following.
	
	Now we choose $\Phi$ and $\lambda(t)$ such that, for any $0<\sigma'<\sigma<1, \ w=\sigma-\sigma',$\\
\\	
(1)$0\le \Phi \le 1$, supp($\Phi) \subset B_{\sigma r}(o)$, $\Phi =1$ in $B_{\sigma 'r}(o)$ and $|\nabla \Phi|\le2(wr)^{-1}$;\\
(2)$0\le \lambda \le 1$, $\lambda=0$ in $(-\infty, s-\sigma r^2)$, $\lambda=1$ in $( s-\sigma' r^2,+\infty)$, and $|\lambda'(t)|\le 2(wr)^{-2}$.\\
\\
Let $I_\sigma=(s-\sigma r^2, s)$ and $I_{\sigma'}=(s-\sigma' r^2, s)$.  For any $t \in I_{\sigma'}$, integrating the above inequality over $(s-\sigma r^2, t)$, we obtain
	
	\begin{equation}\label{4.2}
		\sup\limits_{I_{\sigma'}} \left(\int_{B}\Phi^2 u^2dv\right)  \le C(wr)^{-2}\int_{Q_\sigma}u^2dvdt, 	
	\end{equation}
	and
	\begin{equation}\label{4.3}
		\int_{B\times I_{\sigma'} }    |\nabla (\Phi u)|^2dvdt  \le C(wr)^{-2}\int_{Q_\sigma}u^2dvdt.
	\end{equation}
	
	On the other hand, by the H\"older inequality and the Sobolev inequality in Theorem \ref{Sobolev2}, for some constant $\mu=4n-2$, we have
	
	\begin{equation}
		\begin{aligned}	\label{4.4}
			\int_{B }  g^{2(1+\frac{2}{\mu})}dv&\le \left(\int_{B }  |g|^{\frac{2\mu}{\mu-2}}dv\right)^{\frac{\mu-2}{\mu}}\left(\int_{B }  |g|^2dv\right)^{\frac{2}{\mu}}\\
			&\le \left
			(\int_{B }  |g|^2dv\right)^{\frac{2}{\mu}}\left(E(B)\int_{B } (|\nabla g|^2+r^{-2}g^2)dv\right)
		\end{aligned}
	\end{equation}
	for all $g \in C^\infty(B)$, where $E(B)=\frac{\tilde{c}_5(n)e^{\tilde{c}_6(n)(Kr^2+L(R))}}{\vol (B_{r}(o))^\frac{2}{\mu}}r^2$.
	
	Setting $g=\Phi u$, $\theta=1+\frac{2}{\mu}$, \eqref{4.4} becomes
	\begin{equation}	
		\int_{B }  (\Phi u)^{2\theta}dv\le\left(	\int_{B } (\Phi u)^2dv\right)^{\frac{2}{\mu}}\left(E(B)\int_{B } (|\nabla (\Phi u)|^2+r^{-2}(\Phi u)^2)dv\right).\nonumber
	\end{equation}
	
	Combining \eqref{4.2}, \eqref{4.3} and integrating the above inequality over $(s-\sigma' r^2, s)$, we obtain
	\begin{equation}	
		\int_{s-\sigma' r^2}^{s}\int_{B }  (\Phi u)^{2\theta}dvdt\ge	
		\int_{Q_{\sigma'} }u^{2\theta}dvdt,\nonumber
	\end{equation}

	\begin{equation}	
		\int_{B }  (\Phi u)^2dv\le
		\sup\limits_{I_{\sigma'}}\int_{B }  (\Phi u)^2dv\le C(wr)^{-2}\int_{Q_\sigma}u^2dv, 		\nonumber
	\end{equation}
	and
	\begin{equation}
		\begin{aligned}		
			\int_{s-\sigma' r^2}^{s}\int_{B } \left(|\nabla (\Phi u)|^2+r^{-2}(\Phi u)^2\right)dvdt&=\int_{B\times I_{\sigma'}}|\nabla(\Phi u)|^2dvdt+\int_{s-\sigma' r^2}^{s}\int_{B } r^{-2}(\Phi u)^2dvdt\\
			&\le\int_{B\times I_{\sigma'}}|\nabla(\Phi u)|^2dvdt+\sigma' r^2\times r^{-2}\sup\limits_{I_{\sigma'}}\int_{B } (\Phi u)^2dv\\
			&\le C(wr)^{-2}\int_{Q_\sigma}u^2dvdt,\nonumber	
		\end{aligned}
	\end{equation}
	which implies
	
	\begin{equation}	
		\int_{Q_{\sigma'}}  u^{2\theta}dvdt\le E(B)\left(C(wr)^{-2}\int_{Q_\sigma}u^2dvdt\right)^{\theta}.\nonumber
	\end{equation}
	For any $m\ge 1$, $u^m$ is also a smooth non-negative subsolution of the heat equation. Hence the above inequality indeed implies
	\begin{equation}\label{4.5}
		\int_{Q_{\sigma'}}  u^{2m\theta}dvdt\le E(B)\left(C(wr)^{-2}\int_{Q_\sigma}u^{2m}dvdt\right)^{\theta}
	\end{equation}
	for $m\ge 1$.
	
	Let $w_i=(1-\delta)2^{-i}$, which satisfies $\sum_{1}^{\infty}w_i=1-\delta$. Let $\sigma_0=1$, $\sigma_{i+1}=\sigma_i-w_{i+1}\\=1-\sum_{j=1}^{i+1}w_j$.
	Applying \eqref{4.5} for $m=\theta^i$, $\sigma=\sigma_i$, $\sigma'=\sigma_{i+1}$, we have
	\begin{equation}	
		\int_{Q_{\sigma_{i+1}}}  u^{2\theta^{i+1}}dvdt\le E(B)
		\left[C^{i+1}((1-\delta)r)^{-2} \int_{Q_{\sigma_i}} u^{2\theta^i}dvdt \right]^{\theta},\nonumber
	\end{equation}

	i.e.,
	\begin{equation}
		\left(\int_{Q_{\sigma_{i+1}}}  u^{2\theta^{i+1}}dvdt\right)^{\theta^{-(i+1)}}\le E(B)^{\theta^{-(i+1)}}
		C^{(i+1)\theta^{-i}}((1-\delta)r)^{-2\theta^{-i}}\left(\int_{Q_{\sigma_i}} u^{2\theta^i}dvdt\right)^{\theta^{-i}}.\nonumber
	\end{equation}
	
	Iterating from $i=0$ to $\infty$, we obtain
	
	\begin{equation}
		\sup\limits_{Q_\delta}u^2\le E(B)^{\sum_{i=0}^{\infty}\theta^{-(i+1)}}
		C^{\sum_{i=0}^{\infty}(i+1)\theta^{-i}}((1-\delta)r)^{-2\sum_{i=0}^{\infty}\theta^{-i}}
		\int_{Q } u^2dvdt. 	\nonumber
	\end{equation}
	Therefore

	\begin{equation}
		\sup\limits_{Q_\delta}u^2\le C(n) E(B)^{\frac{\mu}{2}}((1-\delta)r)^{-(\mu+2)}||u||_{2,Q}^2,\nonumber
	\end{equation}
i.e.,

\begin{equation}\label{4.6}
	\sup\limits_{Q_\delta}u^2\le \frac{\tilde{c}_7(n)e^{\tilde{c}_8(n)(Kr^2+L(R))}}{(1-\delta)^{4n}r^2\vol (B_{r}(o))}\int_{Q} u^2dvdt.
\end{equation}
	
	Formula \eqref{4.6} is in fact an $L^2$-mean value inequality. The case $p\ge2$ immediately follows, since for any smooth non-negative subsolution $u$ of the heat equation, $u^{\frac{p}{2}}$, $p\ge2$ is also a smooth non-negative subsolution of the heat equation.
	
	Next, for $0<p<2$, we will apply \eqref{4.6} to prove \eqref{eqqu4.2} by a different iterative argument. Let $\sigma \in (0,1)$ and $\eta=\sigma+(1-\sigma)/4$. Then \eqref{4.6} implies
	
	\begin{equation}
		\sup\limits_{Q_\sigma}u\le F(B)(1-\sigma)^{(-1-\frac{\mu}{2})}\left(\int_{Q_\eta}u^2dvdt\right)^\frac{1}{2},\nonumber
	\end{equation}
	where $F(B)=\frac{\tilde{c}_9(n)e^{\tilde{c}_{10}(n)(Kr^2+L(R))}}{r\vol (B_{r}(o))^\frac{1}{2}}$.
	
	\
	Since $\left(\int_{Q_\eta}u^2dvdt\right)^\frac{1}{2}=\left(\int_{Q_\eta}u^pu^{2-p}dvdt\right)^\frac{1}{2}\le\sup\limits_{Q_\eta}u^{1-\frac{p}{2}} \left(\int_{Q}u^pdvdt\right)^\frac{1}{2}$, we have
	
	\begin{equation}\label{4.7}
		||u||_{\infty,Q_\sigma}	=\sup\limits_{Q_\sigma}u\le F(B)(1-\sigma)^{(-1-\frac{\mu}{2})}||u||_{p,Q}^{\frac{p}{2}}||u||_{\infty,Q_\eta}^{1-\frac{p}{2}}.
	\end{equation}
	Now fix $\delta \in (0,1)$ and let   $\sigma_0=\delta$, $\sigma_{i+1}=\sigma_i+(1-\sigma_i)/4$, which satisfy   $1-\sigma_i=(\frac{3}{4})^i(1-\delta)$. Applying \eqref{4.7} to $\sigma=\sigma_i$, and $\eta=\sigma_{i+1}$, we have
	
	\begin{equation}
		||u||_{\infty,Q_{\sigma _i}}\le (\frac{4}{3})^{(1+\frac{\mu}{2})i}F(B)||u||_{p,Q}^{\frac{p}{2}}(1-\delta)^{(-1-\frac{\mu}{2})}||u||_{\infty,Q_{\sigma _{i+1}}}^{1-\frac{p}{2}}.\nonumber
	\end{equation}
	Therefore, for any $i$,
	
	\begin{equation}
		||u||_{\infty,Q_{\delta}}\le (\frac{4}{3})^{(1+\frac{\mu}{2})\sum_{j} j(1-\frac{p}{2})^j}\times[	F(B)||u||_{p,Q}^{\frac{p}{2}}(1-\delta)^{(-1-\frac{\mu}{2})}]^{\sum_{j} (1-\frac{p}{2})^j}||u||_{\infty,Q_{\sigma _i}}^{(1-\frac{p}{2})^i},\nonumber
	\end{equation}
	where $\sum$ denotes the summations from $0$ to $i-1$. Letting $i\to\infty$ we get
	
	\begin{equation}
		||u||_{\infty,Q_{\delta}}\le \left(\frac{4}{3}\right)^{\tilde{c}(n,p)}\times[	F(B)||u||_{p,Q}^{\frac{p}{2}}(1-\delta)^{(-1-\frac{\mu}{2})}]^\frac{2}{p},\nonumber
	\end{equation}
	that is
	\begin{equation}
		\sup\limits_{Q_\delta}u^p\le \frac{\tilde{c}_{11}(n,p)e^{\tilde{c}_{12}(n)(Kr^2+L(R))}}{(1-\delta)^{4n}r^2\vol( B_{r}(o))} \int_{Q} u^pdvdt,\ 0<p<2.\nonumber
	\end{equation}
	Then the proposition follows.\qed\\
	
\end{proposition}
To get the Gaussian upper bound of the heat kernel, let us first recall Davies' double integral estimate \cite{Da}.

\begin{lemma}\label{Davies}(Davies \cite{Da})
	Let $(M^n,g)$ be a complete Riemannian manifold and $H(x,y,t)$ the heat kernel. Let $\mu_1(M)\ge 0$ be the greatest lower bound for the $L^2$-spectrum of the Laplacian $\Delta$ on $M$. Assume that $B_1$ and $B_2$ are bounded subsets of $M$. Then
	
	\begin{equation}
		\int_{B_1}\int_{B_2}H(x,y,t)dydx\le \vol(B_1)^{\frac{1}{2}}\vol(B_2)^{\frac{1}{2}}e^{\left(-\frac{d^2(B_1,B_2)}{4t}-\mu_1(M)t\right)},\nonumber
	\end{equation}	
	where $d(B_1,B_2)$ denotes the distance between $B_1$ and $B_2$.
\end{lemma}

Now we prove the Gaussian upper bound of the heat kernel by applying the mean value inequality in Proposition \ref{4.1} and Lemma \ref{Davies}.\\

\noindent{\it Proof of Theorem \ref{Gaussian}.}  Fixing $x\in B_{\frac{1}{2}R}(o)$ and applying  Proposition \ref{mean value inequality} to the heat kernel $H(x,y,t)$ with $Q=B_{\sqrt{t}}(y)\times(t-(\sqrt{t})^2,t), \delta=\frac{1}{8}$ and $\delta'=\frac{1}{4}$, we have
\begin{equation}\label{4.12}
	\begin{aligned}
		H(x,y,t)\le \sup\limits_{(z,s)\in Q_\delta}H(x,z,s) &\le\frac{\overline{c}_1(n)e^{\overline{c}_2(n)(Kt+\sup\limits_{B_{3\sqrt{t}}(y)}| f|)}}{\left(\frac{1}{8}\right)^{4n}t\vol( B_{\sqrt{t}}(y))}\int_{t-\frac{1}{4}t}^{t}\int_{B_{\frac{1}{4}\sqrt{t}}(y)}H(x,z,s)dz ds\\
		&\le\frac{\overline{c}_1(n)e^{\overline{c}_2(n)(Kt+L(R))}}{\left(\frac{1}{8}\right)^{4n}t\vol( B_{\sqrt{t}}(y))}\int_{t-\frac{1}{4}t}^{t}\int_{B_{\frac{1}{4}\sqrt{t}}(y)}H(x,z,s)dz ds\\
		& \le \frac{\overline{c}_1(n)e^{\overline{c}_2(n)(Kt+L(R))}}{\left(\frac{1}{8}\right)^{4n}t\vol (B_{\sqrt{t}}(y))}\int_{\frac{3}{4}t}^{t}\int_{B_{\sqrt{t}}(y)}H(x,z,s)dz ds\\
		&=\frac{\overline{c}_3(n)e^{\overline{c}_2(n)(Kt+L(R))}}{\vol (B_{\sqrt{t}}(y))}\int_{B_{\sqrt{t}}(y)}H(x,z,s')dz
	\end{aligned}
\end{equation}
for some $s'\in \left(\frac{3}{4}t,t\right)$, where $Q_{\delta}=B_{\frac{1}{8}\sqrt{t}}(y)\times(t-\frac{1}{8}(\sqrt{t})^2,t)$ and $B_{\sqrt{t}}(y)\subset B_{R}(o)$ for any $y\in B_{\frac{1}{2}R}(o)$ and $0<t<\frac{R^2}{4}$.

Fixing $z\in B_{\sqrt{t}}(y)$ and applying Proposition \ref{mean value inequality} to the heat kernel $H(x,z,s')$ with $Q=B_{\sqrt{t}}(x)\times(t-(\sqrt{t})^2,t)$, $\delta=\frac{1}{4}$ and $\delta'=\frac{1}{2}$, we have
\begin{equation}\label{4.13}
	\begin{aligned}
		H(x,z,s')\le \sup\limits_{(x',s'')\in Q_\delta}H(x',z,s'')&\le\frac{\overline{c}_4(n)e^{\overline{c}_5(n)(Kt+L(R))}}{\left(\frac{1}{4}\right)^{4n}t\vol( B_{\sqrt{t}}(x))}\int_{t-\frac{1}{2}t}^{t}\int_{B_{\frac{1}{2}\sqrt{t}}(x)}H(w,z,\tau)dw d\tau\\
		& \le \frac{\overline{c}_4(n)e^{\overline{c}_5(n)(Kt+L(R))}}{\left(\frac{1}{4}\right)^{4n}t\vol( B_{\sqrt{t}}(x))}\int_{\frac{1}{2}t}^{t}\int_{B_{\sqrt{t}}(x)}H(w,z,\tau)dw d\tau\\
		&=\frac{\overline{c}_6(n)e^{\overline{c}_5(n)(Kt+L(R))}}{\vol (B_{\sqrt{t}}(x))}\int_{B_{\sqrt{t}}(x)}H(w,z,\tau')dw
	\end{aligned}
\end{equation}
for some $\tau'\in\left(\frac{1}{2}t,t\right)$, where $Q_{\delta}=B_{\frac{1}{4}\sqrt{t}}(x)\times(t-\frac{1}{4}(\sqrt{t})^2,t)$ and $B_{\sqrt{t}}(x)\subset B_{R}(o)$ for any $x\in B_{\frac{1}{2}R}(o)$ and $0<t<\frac{R^2}{4}$.\\
Combining \eqref{4.12} and \eqref{4.13}, the heat kernel satisfies
\begin{equation}\label{4.14}
	H(x,y,t)\le \frac{\overline{c}_7(n)e^{\overline{c}_8(n)(Kt+L(R))}}{\vol (B_{\sqrt{t}}(x))\vol (B_{\sqrt{t}}(y))}\int_{B_{\sqrt{t}}(y)}\int_{B_{\sqrt{t}}(x)}H(w,z,\tau')dw dz
\end{equation}
for any $x,\ y\in B_{\frac{1}{2}R}(o)$ and $0<t<\frac{R^2}{4}$.\\
Using Lemma \ref{Davies} and noticing that $\tau'\in \left(\frac{1}{2}t,t\right)$, \eqref{4.14} becomes
\begin{equation}\label{4.15}
		H(x,y,t)\le \frac{\overline{c}_7(n)e^{\overline{c}_8(n)(Kt+L(R))}}{\vol (B_{\sqrt{t}}(x))^{^{\frac{1}{2}}}\vol (B_{\sqrt{t}}(y))^{^{\frac{1}{2}}}}e^{\left(-\frac{d^2(B_{\sqrt{t}}(x),B_{\sqrt{t}}(y))}{4t}\right)}
\end{equation}
for all $x,\ y\in B_{\frac{1}{2}R}(o)$ and $0<t<\frac{R^2}{4}$.\\
Notice that if $d(x,y)\le2\sqrt{t}$, then $d(B_{\sqrt{t}}(x),B_{\sqrt{t}}(y))=0$, and hence
\begin{equation}
	-\frac{d^2(B_{\sqrt{t}}(x),B_{\sqrt{t}}(y))}{4t}=0\le-\frac{d^2(x,y)}{4t},\nonumber
\end{equation}
and if $d(x,y)>2\sqrt{t}$, then $d(B_{\sqrt{t}}(x),B_{\sqrt{t}}(y))=d(x,y)-2\sqrt{t}$, and hence
\begin{equation}
	-\frac{d^2(B_{\sqrt{t}}(x),B_{\sqrt{t}}(y))}{4t}=-\frac{(d(x,y)-2\sqrt{t})^2}{4t}\le1-\frac{d^2(x,y)}{4(1+\epsilon)t}+\frac{1}{\epsilon}\nonumber
\end{equation}
for any $\epsilon>0$.
Combining the above two conditions, gives that
\begin{equation}
e^{\left(-\frac{d^2(B_{\sqrt{t}}(x),B_{\sqrt{t}}(y))}{4t}\right)}\le e^{\left(	-\frac{d^2(x,y)}{4(1+\epsilon)t}+1+\frac{1}{\epsilon}\right)}.\nonumber
\end{equation}
Therefore, by \eqref{4.15} we have
\begin{equation}
H(x,y,t)\le\frac{\overline{c}_9(n,\epsilon)e^{\overline{c}_8(n)(Kt+L(R))}}{\vol (B_{\sqrt{t}}(x))^{\frac{1}{2}}\vol (B_{\sqrt{t}}(y))^{\frac{1}{2}}}e^{\left(-\frac{d^2(x,y)}{(4+\epsilon)t}\right)}
\end{equation}
for all $x,y \in B_{\frac{1}{2}R}(o)$ and $0<t<\frac{R^2}{4}$, where $\lim\limits_{\epsilon \to 0} \overline{c}_9(n,\epsilon)=\infty$. \eqref{4.9} follows by letting $R \to \infty$.\qed \\

Notice that Theorem \ref{volume element comparison} implies
\begin{equation}
	\begin{aligned}
			\vol (B_{\sqrt{t}}(y))&\le \vol (B_{\sqrt{t}+d(x,y)}(x))\\
			&\le\left(\frac{\sqrt{t}+d(x,y)}{\sqrt{t}}\right)^ne^{\frac{K}{6}(d^2(x,y)+2d(x,y)\sqrt{t})+6L(R)}\vol (B_{\sqrt{t}}(x)) \nonumber
	\end{aligned}
\end{equation}
for all $x,y\in B_{\frac{1}{4}R}(o)$ and $B_{\sqrt{t}+d(x,y)}(x) \subset B_{R}(o)$ with $0<t<\frac{R^2}{16}$.
Therefore, the upper bound in Theorem \ref{Gaussian} can be rewritten as follows.

\begin{corollary} \label{cor4.4}
	Under the same assumptions as Theorem \ref{Gaussian}, for any $\epsilon>0$, there exist constants $C_5(n,\epsilon)$ and $C_6(n)$, such that
	\begin{equation}\label{4.11}
		H(x,y,t)\le\frac{C_5(n,\epsilon)e^{C_6(n)(KR^2+L(R)+Kd^2(x,y))}}{\vol(B_{\sqrt{t}}(y)) }\left(1+\frac{d(x,y)}{\sqrt{t}}\right)^{\frac{n}{2}}e^{\left(-\frac{d^2(x,y)}{(4+\epsilon)t}\right)}
	\end{equation}
	for all $x,y \in B_{\frac{1}{4}R}(o)$ and $0<t<R^2/16$, where $L(R)=\sup\limits_{B_{3R}(o)}|f|$ and $\lim\limits_{\epsilon \to 0} C_5(n,\epsilon)=\infty$.
	
\end{corollary}


\section{  $L^1$-Liouville theorem and uniqueness of $L^1$ solutions of the heat equation  }
In this section, inspired by the work of Li \cite{Li}, we prove the $L^1$-Liouville theorem for non-negative subharmonic functions and the uniqueness theorem for $L^1$-solutions of the heat equation on complete noncompact Riemannian manifolds with Bakry-\'Emery Ricci curvature bounded below and the potential function of at most quadratic growth.


We start from a useful lemma.
\begin{lemma}(Theorem 11.8 in \cite{Gr})\label{stochastically complete lemma}
	Let $(M^n,g)$ be a complete Riemannian manifold. If, for some point $x_0\in M$,
	\begin{equation}
		\int_{1}^{+\infty}\frac{R\,dR}{\ln(\vol (B_{R}(x_0)))}=\infty,\nonumber
	\end{equation}	
	then	$(M^n,g)$ is stochastically complete, i.e.,
		\begin{equation}
		\int_{M}	H(x,y,t)dy=1.\nonumber
	\end{equation}
\end{lemma}

Under our assumptions mentioned earlier, it is easy to check the stochastic completeness of the manifold due to Theorem \ref{volume element comparison}.

\begin{proposition}\label{Gri}
		Let $(M^n,g)$ be a complete non-compact Riemannian manifold with $\Ric_f\ge -Kg$ for some constant $K\ge0$. Assume  there exist non-negative constants $a$ and $b$ such that
	\begin{equation}	
		|f|(x)\le ar^2(x)+b \ for \ all \ x\in M,\nonumber
	\end{equation}	
	where $r(x)=d(x,o)$ is the geodesic distance function to a fixed point $o\in M$. Then	$(M^n,g)$ is stochastically complete.

\proof	In \eqref{VC},\ letting $r_1\to 0, \ r_2=R$, and $p=o\in M$ yields
	\begin{equation}
	\vol (B_{R}(o))\le c(n,b)R^ne^{c(K,a)R^2}\nonumber
\end{equation}	
for all $R>1$. Hence	
	\begin{equation}
		\int_{1}^{+\infty}\frac{R\,dR}{\ln(\vol (B_{R}(o)))}=\infty.\nonumber
\end{equation}	
By Lemma \ref{stochastically complete lemma}, this implies that $(M^n, g)$ is stochastically complete.\qed\\	
\end{proposition}

Now, we are ready to check the integration by parts formula by using the upper bound of the heat kernel in Corollary \ref{cor4.4} and the mean value inequality in Proposition \ref{4.1}.
\begin{proposition}\label{integration by parts formula}
	Under the same assumptions as Proposition \ref{Gri}, for any non-negative $L^1$-integrable subharmonic function $u$, we have	
\begin{equation}
	\int_{M} \Delta_y	H(x,y,t)u(y)dy=\int_{M} H(x,y,t)\Delta_yu(y)dy\nonumber
\end{equation}	
for any $x\in M$ and $t>0$.

\proof Applying integration by parts on $ B_R(o)$  gives
\begin{equation}\label{Green's formula}
	\begin{aligned}
	&\left|\int_{B_R(o)} \Delta_y	H(x,y,t)u(y)dy-\int_{B_R(o)} H(x,y,t)\Delta_yu(y)dy\right| \\
	&=\left|\int_{\partial B_R(o)} \frac{\partial}{\partial r}	H(x,y,t)u(y)dS-\int_{\partial B_R(o)} H(x,y,t)\frac{\partial}{\partial r}u(y)dS\right|\\
	&\le \int_{\partial B_R(o)}  |\nabla_y H|(x,y,t)u(y)dS+\int_{\partial B_R(o)}  H(x,y,t)|\nabla_y u|(y)dS   ,
		\end{aligned}
\end{equation}
where $dS$ denotes the area measure on $\partial B_R(o)$. We shall show that the above two boundary integrals vanish as $R\to \infty$.
Without loss of generality, we assume $R>1$ and $x\in B_{\frac{1}{8}R}(o)$.

Step 1.
Since $	|f|(x)\le ar^2(x)+b $, by Proposition \ref{mean value inequality}, we get
	\begin{equation}\label{5.1}
	\sup\limits_{B_R(o)}u\le \frac{\tilde{C}_1(n)e^{\tilde{C}_2(n,K,a,b)R^2}}{\vol (B_{2R}(o))}\int_{B_{2R}(o)} u\le \frac{\tilde{C}_1(n)e^{\tilde{C}_2(n,K,a,b)R^2}}{\vol (B_{2R}(o))}||u||_1.
\end{equation}
Let $\phi(y)=\phi(d(y,o))$ be a non-negative cut-off function satisfying
$0\le \phi\le1$, \ $|\nabla \phi|\le \sqrt{3}$,  $\phi(y)=1$  on $B_{R+1}(o)\backslash B_R(o)$, and $\phi(y)=0$ on $B_{R-1}(o)\cup (M\backslash B_{R+2}(o))$.
Since $u$ is subharmonic, by the Cauchy-Schwarz inequality, we have
	\begin{equation}
		\begin{aligned}
		0\le \int_{M} \phi^2u\Delta u&=-\int_{M} \left<\nabla (\phi^2u),\nabla u\right>\\
		&=-2\int_{M} \phi u\left<\nabla \phi,\nabla u\right>-\int_{M}\phi^2|\nabla u|^2\\
		&\le 2\int_{M} |\nabla \phi|^2u^2-\frac{1}{2}\int_{M}\phi^2|\nabla u|^2,\nonumber
	\end{aligned}
	\end{equation}
i.e.,
$$
\int_{M}\phi^2|\nabla u|^2 \le 4\int_{M} |\nabla \phi|^2u^2.
$$

It then follows from \eqref{5.1} that
	\begin{equation}
	\begin{aligned}
	\int_{B_{R+1}(o)\backslash B_R(o)} |\nabla u|^2\le\int_{M}\phi^2|\nabla u|^2 &\le 4\int_{M} |\nabla \phi|^2u^2\\
	&\le 12\int_{B_{R+2}(o)} u^2\\
	&\le 12\sup\limits_{B_{R+2}(o)}u \times ||u||_{1}\\
	&\le \frac{\tilde{C}_3(n)e^{\tilde{C}_4(n,K,a,b)(R+2)^2}}{\vol (B_{2R+4}(o))}||u||_1^2.\nonumber
	\end{aligned}
\end{equation}
On the other hand, the Cauchy-Schwarz inequality implies that
	\begin{equation}
	\int_{B_{R+1}(o)\backslash B_R(o)} |\nabla u|\le\left(	\int_{B_{R+1}(o)\backslash B_R(o)} |\nabla u|^2\right)^{\frac{1}{2}}\cdot
\left[\vol (B_{R+1}(o))-\vol (B_R(o))	\right]^{\frac{1}{2}}.\nonumber	
\end{equation}
Combining the above two inequalities, we have
\begin{equation}\label{5.2}
		\int_{B_{R+1}(o)\backslash B_R(o)} |\nabla u|\le
		\tilde{C}_5(n)e^{\tilde{C}_6(n,K,a,b)(R+2)^2}||u||_1.
\end{equation}

Step 2. By letting $\epsilon=1$ in Corollary \ref{cor4.4}, the heat kernel $H(x,y,t)$ satisfies
\begin{equation}\label{5.3}
	H(x,y,t)\le\frac{\tilde{C}_{7}(n)e^{\tilde{C}_{8}(n,K,a,b)[R^2+d^2(x,y)]}}{\vol  (B_{\sqrt{t}}(x))}\left(1+\frac{d(x,y)}{\sqrt{t}}\right)^{\frac{n}{2}}e^{\left(-\frac{d^2(x,y)}{5t}\right)}
\end{equation}
for all $x,y \in B_{\frac{1}{4}R}(o)$ and $0<t<\frac{R^2}{16}$.\\
Together with \eqref{5.2} we get
\begin{equation}
	\begin{aligned}
		J_1:&= \int_{B_{R+1}(o)\backslash B_R(o)} H(x,y,t)|\nabla u|(y)dy \\
		&\le \sup\limits_{y\in B_{R+1}(o)\backslash B_R(o)}H(x,y,t)\cdot\int_{B_{R+1}(o)\backslash B_R(o)}|\nabla u|(y)dy\\
		&\le\frac{\tilde{C}_{9}(n)e^{\tilde{C}_{10}(n,K,a,b)[(R+2)^2+(d(x,o)+R+1)^2]-\frac{(R-d(x,o))^2}{5t}}||u||_1}{\vol (B_{\sqrt{t}}(x))}\left(1+\frac{d(x,o)+R+1}{\sqrt{t}}\right)^{\frac{n}{2}}.
		\nonumber
	\end{aligned}
\end{equation}
Thus, for $T$ sufficiently small, all $t\in(0,T)$, and $d(x,o)\le \frac{R}{8}$, there exists a fixed constant $\beta>0$ such that
\begin{equation}
	\begin{aligned}
		J_1\le\frac{\hat{C}_1(n,K,a,b)||u||_1}{\vol (B_{\sqrt{t}}(x))}\left(1+\frac{2R+1}{\sqrt{t}}\right)^{\frac{n}{2}}e^{\big[-\beta R^2+\hat{C}_2(n,K,a,b)R\big]}.\nonumber
	\end{aligned}
\end{equation}
Hence for all $t\in(0,T)$ and all $x\in M$, $J_1\rightarrow 0$ as $R\rightarrow\infty$.

   Step 3. We show that $\displaystyle \int_{B_{R+1}(o)\backslash B_{R}(o)} |\nabla_y H|(x,y,t)u(y)dy\to 0$ as $R\to \infty$. Frist, consider the integral
\begin{equation}
	\begin{aligned}
		\int_{M} \phi^2(y)|\nabla_y H|^2(x,y,t)dy&=-2\int_{M}\left<H(x,y,t)\nabla \phi(y),\phi(y)\nabla_y H(x,y,t)\right>dy\\
		&\ \ \ \ -\int_{M}\phi^2(y)H(x,y,t)\Delta_yH(x,y,t)dy\\
		&\le2\int_{M}|\nabla\phi|^2(y)H^2(x,y,t)dy +\frac{1}{2}\int_{M}\phi^2(y)|\nabla_y H|^2(x,y,t)dy\\
		&\ \ \ -\int_{M}\phi^2(y)H(x,y,t)\Delta_yH(x,y,t)dy.\nonumber
	\end{aligned}
\end{equation}
This implies
\begin{equation}\label{5.4}
	\begin{aligned}
		&\int_{B_{R+1}(o)\backslash B_R(o)}|\nabla_y H|^2(x,y,t)dy \le \int_{M}\phi^2(y)|\nabla_y H|^2(x,y,t)dy \\
		&\le4\int_{M}|\nabla \phi|^2(y)H^2(x,y,t)dy -2\int_{M}\phi^2(y)H(x,y,t)\Delta_yH(x,y,t)dy\\
		&\le12\int_{B_{R+2}(o)\backslash B_{R-1}(o)}H^2(x,y,t)dy+2\int_{B_{R+2}(o)\backslash B_{R-1}(o)}H(x,y,t)|\Delta_yH|(x,y,t)dy\\
		&\le12\int_{B_{R+2}(o)\backslash B_{R-1}(o)}H^2(x,y,t)dy+2\left(\int_{B_{R+2}(o)\backslash B_{R-1}(o)}H^2\right)^{\frac{1}{2}}\left(\int_{M}|\Delta_yH|^2(x,y,t)\right)^{\frac{1}{2}}.
	\end{aligned}
\end{equation}
By Proposition \ref{Gri}, we have
\begin{equation}
	\int_MH(x,y,t)dy=1 \nonumber
\end{equation}
for all $x\in M$ and $t>0$. Combining this with \eqref{5.3} gives
\begin{equation}\label{5.5}
	\begin{aligned}
		&\int_{B_{R+2}(o)\backslash B_{R-1}(o)}H^2(x,y,t)dy\le\sup\limits_{y\in B_{R+2}(o)\backslash B_{R-1}(o)}H(x,y,t)\\
		&\le\frac{\tilde{C}_{11}(n)e^{\tilde{C}_{12}(n,K,a,b)[(R+2)^2+(d(x,o)+R+2)^2]-\frac{(R-1-d(o,x))^2}{5t}}}{\vol (B_{\sqrt{t}}(x))}\left(1+\frac{d(x,o)+R+2}{\sqrt{t}}\right)^{\frac{n}{2}}.
	\end{aligned}
\end{equation}
Form (29.15) in \cite{Li2}, we know that
\begin{equation}\label{claim}
	\int_{M}(\Delta_yH)^2(x,y,t)dy \le \frac{\tilde{C}}{t^2}H(x,x,t).
\end{equation}
Combining \eqref{5.4}, \eqref{5.5} and \eqref{claim}, we obtain
\begin{equation}
	\begin{aligned}
		&\int_{B_{R+1}(o)\backslash B_R(o)}|\nabla_y H|^2(x,y,t)dy \le \tilde{C}_{13}(n)e^{\tilde{C}_{14}(n,K,a,b)[(R+2)^2+(d(x,o)+R+2)^2]-\frac{(R-1-d(x,o))^2}{10t}}\\
		&\times\left[\vol (B_{\sqrt{t}}(x))^{-1}+\vol (B_{\sqrt{t}}(x))^{-\frac{1}{2}}t^{-1}H^{\frac{1}{2}}(x,x,t)\right]\left(1+\frac{d(x,o)+R+2}{\sqrt{t}}\right)^{\frac{n}{2}}.\nonumber
	\end{aligned}
\end{equation}
By the Cauchy-Schwarz inequality, we get
\begin{equation}\label{5.7}
	\begin{aligned}
		&\int_{B_{R+1}(o)\backslash B_{R}(o)}|\nabla_y H|(x,y,t)dy\\
		&\le\left[\vol (B_{R+1}(o))-\vol (B_{R}(o))\right]^{\frac{1}{2}}\cdot\left(\int_{B_{R+1}(o)\backslash B_{R}(o)}|\nabla_y H|^2(x,y,t)dy\right)^{\frac{1}{2}}\\
		&\le \vol (B_{R+1}(o))^{\frac{1}{2}}\tilde{C}_{16}(n)e^{\tilde{C}_{15}(n,K,a,b)[(R+2)^2+(d(x,o)+R+2)^2]-\frac{(R-1-d(x,o))^2}{20t}}\\
		&\times \left[\vol (B_{\sqrt{t}}(x))^{-1}+\vol  (B_{\sqrt{t}}(x))^{-\frac{1}{2}}t^{-1}H^{\frac{1}{2}}(x,x,t)\right]^{\frac{1}{2}} \left(1+\frac{d(x,o)+R+2}{\sqrt{t}}\right)^{\frac{n}{4}}.
	\end{aligned}
\end{equation}
Therefore, by \eqref{5.1} and \eqref{5.7}, we have
\begin{equation}
	\begin{aligned}
	J_2: &=	\int_{B_{R+1}(o)\backslash B_R(o)}|\nabla_y H|(x,y,t)u(y)dy\\
		&\le \sup\limits_{y \in B_{R+1}(o)\backslash B_R(o)}u(y)\cdot \int_{B_{R+1}(o)\backslash B_R(o)}|\nabla_y H|(x,y,t)dy\\
			&\le \frac{\tilde{C}_{18}(n)e^{\tilde{C}_{17}(n,K,a,b)(R+1)^2}}{\vol (B_{2R+2}(o))^{\frac{1}{2}}}||u||_1
			e^{\tilde{C}_{15}(n,K,a,b)[(R+2)^2+(d(x,o)+R+2)^2]-\frac{(R-1-d(x,o))^2}{20t}}\\
			&\times\left[\vol (B_{\sqrt{t}}(x))^{-1}+\vol  (B_{\sqrt{t}}(x))^{-\frac{1}{2}}t^{-1}H^{\frac{1}{2}}(x,x,t)\right]^{\frac{1}{2}} \left(1+\frac{d(x,o)+R+2}{\sqrt{t}}\right)^{\frac{n}{4}}\\
			&\le \frac{\tilde{C}_{18}(n)e^{\tilde{C}_{19}(n,K,a,b)[(R+2)^2+(d(x,o)+R+2)^2]-\frac{(R-1-d(x,o))^2}{20t}}||u||_1}{\vol (B_{2}(o))^{\frac{1}{2}}}\\
			&\times\left[\vol (B_{\sqrt{t}}(x))^{-1}+\vol  (B_{\sqrt{t}}(x))^{-\frac{1}{2}}t^{-1}H^{\frac{1}{2}}(x,x,t)\right]^{\frac{1}{2}} \left(1+\frac{d(x,o)+R+2}{\sqrt{t}}\right)^{\frac{n}{4}}.\nonumber
	\end{aligned}
\end{equation}

Similar to the case of $J_1$, by  choosing $T$ sufficiently small, for all $t\in (0, T)$ and all $x\in M$, $J_2$ also tends to zero when $R$ tends to infinity.

Step 4. By the mean value theorem, for all $R>0$ there exists $\bar{R}\in(R,R+1)$ such that
\begin{equation}
	\begin{aligned}
J:&=\int_{\partial B_{\bar{R}}(o)}  |\nabla_y H|(x,y,t)u(y)dS+\int_{\partial B_{\bar{R}}(o)}  H(x,y,t)|\nabla_y u|(y)dS \\
&=	\int_{B_{R+1}(o)\backslash B_{R}(o)} |\nabla_y H|(x,y,t)u(y)dy+\int_{B_{R+1}(o)\backslash B_{R}(o)}  H(x,y,t)|\nabla_y u|(y)dy \\
&=J_2+J_1.\nonumber
 	\end{aligned}
\end{equation}

From step 2 and step 3, we know that by choosing $T$ sufficiently small, for all $t\in (0, T)$ and all $x\in M$, $J$  tends to zero along a sequence of radii tending to infinity. Since $H(x,y,t)\Delta u\geq 0$ and $\Delta_y H(x,y,t)=\frac{\partial }{\partial t}H(x,y,t)$ is bounded by a result of Grigor'yan \cite{Gr2}, one can see from \eqref{Green's formula} that $J\rightarrow 0$ for all $R\rightarrow \infty$. Therefore, Proposition \ref{integration by parts formula} is true for $t$ sufficiently small.

Step 5. For all $t\in(0, T)$ and $s\in (0,+\infty)$, using the semigroup property of the heat kernel, we have
\begin{equation}
	\begin{aligned}
		\int_{M} \Delta_y	H(x,y,t+s)u(y)dy
		 &=\int_{M}\int_{M} H(x,z,s)\Delta_yH(z,y,t)dzu(y)dy\\
		  &=\int_{M}\left(\int_{M}\Delta_yH(z,y,t)u(y)dy\right)H(x,z,s)dz\\
		   &=\int_{M}\left(\int_{M}H(z,y,t)\Delta_yu(y)dy\right)H(x,z,s)dz\\
		     &=\int_{M}	H(x,y,t+s)\Delta_yu(y)dy .\nonumber	
	\end{aligned}
\end{equation}
This completes the proof of Proposition \ref{integration by parts formula} for all time $t>0$.\qed\\
\end{proposition}

Applying the regularity theory of harmonic functions, combining Proposition \ref{Gri} and Propositon \ref{integration by parts formula}, we can obtain the $L^1$-Liouville property.\\

\noindent{\it Proof of Theorem \ref{Liouville theorem}.} Let $u(x)$ be a non-negative $L^1$-integrable subharmonic function on $M$. We define a space-time function
\begin{equation}
	u(x,t)=\int_{M}H(x,y,t)u(y)dy\nonumber
\end{equation}	
with initial data $	u(x,0)=u(x)$. From Proposition \ref{integration by parts formula}, we conclude that
\begin{equation}
	\begin{aligned}
	\frac{\partial}{\partial t}u(x,t)&=\int_{M}\frac{\partial}{\partial t}H(x,y,t)u(y)dy\\
	&=\int_{M}\Delta_yH(x,y,t)u(y)dy\\
	&=\int_{M}H(x,y,t)\Delta_yu(y)dy\ge0,\nonumber
		\end{aligned}
\end{equation}
that is, $u(x,t)$ is increasing in $t$. By Proposition \ref{Gri},
\begin{equation}
	\int_{M}H(x,y,t)dy=1\nonumber
\end{equation}	
for all $x\in M$ and $t>0$. So we have
\begin{equation}
	\int_{M}u(x,t)dx=	\int_{M}\int_{M}H(x,y,t)u(y)dydx=	\int_{M}u(y)dy.\nonumber
\end{equation}	
Since $u(x,t)$ is increasing in $t$, so $u(x,t)=u(x)$ and hence $u(x)$ is a non-negative harmonic function, i.e., $\Delta u(x)=0$.

On the other hand, for any positive constant $\epsilon$, let us define a new function $h(x)=\min\{u(x),\epsilon\}$. Then $h$ satisfies

  \ $0\le h(x)\le u(x)$, \ $|\nabla h|\le |\nabla u|$ \ and $\Delta h(x)\le 0$.

 So $h$ has the integration
 by parts formula as $u$. Similarly we define $h(x,t)=\int_{M}H(x,y,t)h(y)dy$ and
\begin{equation}
	\begin{aligned}
		\frac{\partial}{\partial t}h(x,t)&=\int_{M}\frac{\partial}{\partial t}H(x,y,t)h(y)dy\\
		&=\int_{M}H(x,y,t)\Delta_yh(y)dy\le 0.\nonumber
	\end{aligned}
\end{equation}
By the same argument, we have that $\Delta h(x)=0$.

By the regularity theory of harmonic functions, that is impossible unless $h=u$ or $h=\epsilon$. Since $\epsilon$ is arbitrary and $u$ is non-negative, so $u$ must be identically constant. The theorem then follows from the fact that the absolute value of a harmonic function is a non-negative subharmonic function.\qed\\

With the $L^1$-Liouville property, one can prove the uniqueness of $L^1$ solution of the heat equation. \\

\noindent{\it Proof of Theorem \ref{L^1}.} Let $u(x,t)\in L^1$ be a non-negative function satisfying the assumptions in Theorem \ref{L^1}. For $\epsilon>0$, we define a space-time function
\begin{equation}\label{5.9}
	u_\epsilon(x,t)=\int_MH(x,y,t)u(y,\epsilon)dy
\end{equation}
and
\begin{equation}
	F_\epsilon(x,t)=\max\{0,u(x,t+\epsilon)-u_\epsilon(x,t)\}.\nonumber
\end{equation}
Then $F_\epsilon(x,t)$ is non-negative and satisfies
\begin{equation}
	\lim\limits_{t\rightarrow 0}F_\epsilon(x,t)=0,\ \left(\Delta-\frac{\partial}{\partial t}\right)F_\epsilon(x,t)\ge0. \nonumber
\end{equation}
Let $T>0$ be fixed. Let $h(x)=\int_{0}^{T}F_\epsilon(x,t)dt$, which implies
\begin{equation}\label{5.10}
	\Delta h(x)=\int_{0}^{T}\Delta F_\epsilon(x,t)dt\ge\int_{0}^{T}\partial_tF_\epsilon(x,t)dt=F_\epsilon(x,T)\ge0,
\end{equation}
and
\begin{equation}
	\begin{aligned}
		\int_M h(x)dx&=\int_{0}^{T}\int_MF_\epsilon(x,t)dxdt\le \int_{0}^{T}\int_M|u(x,t+\epsilon)-u_\epsilon(x,t) |dxdt \\
		&\le \int_{0}^{T}\int_M u(x,t+\epsilon)dxdt+\int_{0}^{T}\int_M u_\epsilon(x,t) dxdt<\infty,\nonumber
	\end{aligned}	
\end{equation}
where the first term on the right hand is finite from our assumption, and the second term is finite because the heat semigroup is contractive on $L^1$. Therefore, $h(x)$ is a non-negative $L^1$-integrable subharmonic function. By Theorem \ref{Liouville theorem}, $h(x)$ must be constant. Combining with \eqref{5.10} we have $F_\epsilon(x,t)=0$ for all $x\in M$ and $t>0$, which implies
\begin{equation}\label{5.11}
	u_\epsilon(x,t)\ge u(x,t+\epsilon).
\end{equation}
Next we estimate the function $u_\epsilon(x,t) $ in \eqref{5.9}. Applying the upper bound estimate \eqref{4.11} of the heat kernel $H(x,y,t)$ and letting $\epsilon=1$, $R=2d(x,y)+1$, we have
\begin{equation}
	u_\epsilon(x,t)\le\frac{c}{\vol(B_{\sqrt{t}}(x))}\int_M\left[e^{\tilde{c}d^2(x,y)-\frac{d^2(x,y)}{5t}}\left(1+\frac{d(x,y)}{\sqrt{t}}\right)^\frac{n}{2}\right]u(y,\epsilon)dy.\nonumber
\end{equation}
For sufficiently small values of $t>0$, the right-hand side can be estimated by
\begin{equation}
	\frac{C}{\vol(B_{\sqrt{t}}(x))}\int_Mu(y,\epsilon)dy.\nonumber
\end{equation}
Hence as $\epsilon\rightarrow 0$, $u_\epsilon(x,t)\rightarrow 0$ since $\int_Mu(y,\epsilon)dy\rightarrow 0$. However, by the semigroup property, $u_\epsilon(x,t)\rightarrow 0$ for all $x\in M$ and $t>0$. Combining with \eqref{5.11} we get $u(x,t)\le 0$. Therefore $u(x,t)\equiv0$.

To prove that any $L^1$-solution of the heat equation is uniquely determined by its initial data in $L^1$, we suppose that $u_1(x,t),\ u_2(x,t)$ are two $L^1$-integrable solutions of the heat equation $(\Delta-\partial_t)u=0$ with the initial data $u(x,0)\in L^1$. Applying this above result to $v(x,t)=|u_1(x,t)-u_2(x,t)|$, we see that $v(x,t)\equiv 0$. The proof of Theorem \ref{L^1} is finished. \\ \qed

\section{ An $L^{\infty}$ Liouville Property for harmonic functions with polynomial growth }
In this section, we take a detour to prove an $L^{\infty}$ Liouville theorem for harmonic functions with polynomial growth when the Bakry-\'Emery Ricci curvature is non-negative and the potential function is bounded by using a similar idea in \cite{CG}. This is not the main result of the paper, but can be compared with the results in the last section and may be of independent interest.

First, we give a gradient estimate in the integral sense.

\begin{lemma}\label{estimate}
	Let $(M^n,g)$ be a complete non-compact Riemannian manifold with $\Ric_f\ge0$ and $|f|\le L$ on $M$ for some constant $L\ge0$. For any point $o\in M$ and $R>0$, let $u$ be a harmonic function on $B_{2R}(o)$, and $diam\ \partial B_r(o):=\sup\limits_{x,y\in\partial B_r(o)}d(x,y)\le \epsilon r $ with $\epsilon\in(0,\frac{1}{12})$ for all $r\in[R,2R]$. Then we have

\begin{equation}\label{7.1}
	\int_{B_R(o)}|\nabla u|^2 \le \delta^{\frac{1}{\epsilon}-6}\int_{B_{2R}(o)}|\nabla u|^2 ,
\end{equation}
where $\delta=\left(\frac{9c_1e^{c_2L}}{1+9c_1e^{c_2L}}\right)^{\frac{1}{6}}$ and $c_1$, $c_2$ are constants depending only on $n$ from Lemma \ref{Poincare }.
\end{lemma}

\proof Choosing $r\in[R+3\epsilon R,2R-3\epsilon R]$, where $\epsilon\in (0,\frac{1}{12})$, the diameter hypothesis in the Lemma implies that there exists some $x\in \partial B_r(o)$ such that
\begin{equation}
	B_{r+\epsilon R}(o) \backslash B_r(o) \subset B_{\epsilon R+\epsilon r}(x).\nonumber
\end{equation}
Let $\eta (y)$ be a cut-off function with support in $ B_{r+\epsilon R}(o)$ such that
\begin{equation}
	\eta(y)=
	\begin{cases}
		1& \text{ $ y\in B_r(o)$ }, \\
		\frac{r+\epsilon R-d(y,o)}{\epsilon R}& \text{ $ y \in B_{r+\epsilon R}(o) \backslash B_r(o) $ }, \\
		0& \text{ $ y\in M \backslash B_{r+\epsilon R}(o)$ }.\nonumber
	\end{cases}
\end{equation}
We easily observe that
\begin{equation}
	\int_M |\nabla(\eta(u-c))|^2=\int_M \eta^2|\nabla(u-c)|^2+2\eta(u-c)\left<\nabla\eta,\nabla(u-c)\right>+(u-c)^2|\nabla\eta|^2,\nonumber
\end{equation}
where $c$ is a real number to be chosen later.  Notice that
\begin{equation}
	\begin{aligned}
		\int_M \eta^2|\nabla(u-c)|^2+2\eta(u-c)\left<\nabla\eta,\nabla(u-c)\right>&=\int_M\left<\nabla((u-c)\eta^2),\nabla(u-c)\right>\\
		&=-\int_M(u-c)\eta^2\Delta(u-c)=0.\nonumber
	\end{aligned}
\end{equation}
Therefore,
\begin{equation}
	\int_M |\nabla(\eta(u-c))|^2=\int_{B_{r+\epsilon R}(o)}(u-c)^2|\nabla\eta|^2.\nonumber
\end{equation}
According to the definition of $\eta(y)$, by the above equality, we get that
\begin{equation}
	\begin{aligned}
		\int_{B_r(o)}|\nabla u|^2 &\le \int_{B_{r+\epsilon R}(o)}|\nabla(\eta(u-c))|^2 \\
		&=\int_{B_{r+\epsilon R}(o)}(u-c)^2|\nabla\eta|^2\\
		&\le \frac{1}{\epsilon^2R^2}\int_{B_{r+\epsilon R}(o) \backslash B_r(o)}(u-c)^2\\
		&\le \frac{1}{\epsilon^2R^2}\int_{B_{\epsilon R+\epsilon r}(x)}(u-c)^2.\nonumber
	\end{aligned}
\end{equation}
For the right-hand side of the above inequality, if we choose $c=u_{B_{\epsilon R+\epsilon r}(x)}$, then using the Poincar\'e inequality in Lemma \ref{Poincare } and the fact that $r+R\le3R$, we have
\begin{equation}
	\int_{B_r(o)}|\nabla u|^2 \le 9c_1e^{c_2L}\int_{B_{3\epsilon R}(x)}|\nabla u|^2.\nonumber
\end{equation}
This implies that
\begin{equation}
	\int_{B_{r-3\epsilon R}(o)}|\nabla u|^2\le	\int_{B_r(o)}|\nabla u|^2 \le 9c_1e^{c_2L}\int_{B_{r+3\epsilon R}(o)\backslash B_{r-3\epsilon R}(o)}|\nabla u|^2,\nonumber
\end{equation}
that is,
\begin{equation}
	\int_{B_{r-3\epsilon R}(o)}|\nabla u|^2 \le\frac{9c_1e^{c_2L}}{1+9c_1e^{c_2L}}\int_{B_{r+3\epsilon R}(o)}|\nabla u|^2,\ r\in[R+3\epsilon R,2R-3\epsilon R].\nonumber
\end{equation}
Set $r=R+3\epsilon R$, then
\begin{equation}
	\int_{B_{R}(o)}|\nabla u|^2 \le\frac{9c_1e^{c_2L}}{1+9c_1e^{c_2L}}\int_{B_{R+6\epsilon R}(o)}|\nabla u|^2.\nonumber
\end{equation}
Iterating this inequality $N$ times, we finally get
\begin{equation}
	\int_{B_{R}(o)}|\nabla u|^2 \le\left(\frac{9c_1e^{c_2L}}{1+9c_1e^{c_2L}}\right)^N\int_{B_{2R}(o)}|\nabla u|^2\nonumber
\end{equation}
provided that $N6\epsilon R\le R$.
Thus, we can choose $N=\left[\frac{1}{6\epsilon}\right] \ge \frac{1}{6\epsilon}-1>0$ for $\epsilon \in \left(0,\frac{1}{12}\right)$.
\\
Then we have
\begin{equation}
	\int_{B_{R}(o)}|\nabla u|^2 \le\left(\frac{9c_1e^{c_2L}}{1+9c_1e^{c_2L}}\right)^{\frac{1}{6\epsilon}-1}\int_{B_{2R}(o)}|\nabla u|^2.\nonumber
\end{equation}
The desired result follows by choosing $\delta=\left(\frac{9c_1e^{c_2L}}{1+9c_1e^{c_2L}}\right)^{\frac{1}{6}}$. \qed\\

Now we are ready to prove the Liouville property by using Lemma \ref{estimate}.
\begin{theorem}
	Let $(M^n,g)$ be a complete noncompact Riemannian manifold with $\Ric_f\ge0$ and $|f|\le L$. For a base point $o\in M$, if the diameter of the geodesic sphere $\partial B_R(o)$ has a sublinear growth, i.e., $diam\ \partial B_R(o)=o(R),\ R\rightarrow\infty$, then any harmonic function with polynomial growth is constant.
\end{theorem}
\proof Let $u:M \rightarrow \mathbb{R}$ be a harmonic function with polynomial growth of order $\nu$, namely
\begin{equation}
	|u(x)|\le c(1+d(x,o))^\nu.\nonumber
\end{equation}
For $R>>1$, we define
\begin{equation}
	I_R:=\int_{B_R(o)}|\nabla u|^2\nonumber
\end{equation}
and
\begin{equation}
	\epsilon(r):=\sup\limits_{t\ge r}\frac{\rho(t)}{t},\nonumber
\end{equation}
where $\rho(t):=\sup\limits_{x,y\in\partial B_t(o)}d(x,y).$

To estimate $I_R$, we introduce a cut-off function $\xi(x)$ such that
\begin{equation}
	\xi(x)=
	\begin{cases}
		1& \text{ $ x\in B_R(o)$ }, \\
		\frac{2R-d(x,o)}{R}& \text{ $ x \in B_{2R}(o) \backslash B_R(o) $ }, \\
		0& \text{ $ x\in M \backslash B_{2R}(o)$ }.\nonumber
	\end{cases}
\end{equation}
Then we have
\begin{equation}
	\begin{aligned}
		I_R&\le \int_{B_{2R}(o)}|\nabla(\xi u)|^2=\int_{B_{2R}(o)}|u|^2|\nabla \xi|^2+\xi^2|\nabla u|^2+2u\xi\left<\nabla u,\nabla\xi\right>\\
		&=\int_{B_{2R}(o)}|u|^2|\nabla \xi|^2 \le \int_{B_{2R}(o)}c(1+2R)^{2\nu}\frac{1}{R^2}\le \int_{B_{2R}(o)}c(3R)^{2\nu}\frac{1}{R^2}\\
		&=cR^{2\nu-2}\vol(B_{2R}(o)).\nonumber
	\end{aligned}
\end{equation}
The volume comparison Theorem \ref{volume element comparison} shows that
\begin{equation}
	\frac{\vol(B_{2R}(o))}{\vol(B_{1}(o))}\le e^{6L}(2R)^n.\nonumber
\end{equation}
Hence
\begin{equation}\label{7.2}
	I_R\le C(n,L,\vol(B_{1}(o)))R^{2\nu+n-2}.
\end{equation}

On the other hand, if we iterate the inequality \eqref{7.1} proved in Lemma \ref{estimate} $l$ times, we can show that, for all sufficiently large $R$ such that $\epsilon(R)<\frac{1}{12}$,
\begin{equation}
	I_R\le \delta^{-6l+\sum_{j=0}^{l-1}\frac{1}{\epsilon(2^jR)}}I_{2^{l}R},\nonumber
\end{equation}
where $\delta=\left(\frac{9c_1e^{c_2L}}{1+9c_1e^{c_2L}}\right)^\frac{1}{6}$.
Applying \eqref{7.2} to the right-hand side of the above inequality yields
\begin{equation}\label{equuu3.3}
	I_R\le C(n,L,\vol(B_{1}(o)))e^{l\left[\left(\frac{\sum_{j=0}^{l-1}\frac{1}{\epsilon(2^jR)}}{l}-6\right)\ln\delta+(2\nu+n-2)\left(\ln2+\frac{\ln R}{l}\right)\right]}.
\end{equation}
For all sufficiently large $R$,
\begin{equation}
	\lim\limits_{l\rightarrow \infty}\frac{\sum_{j=0}^{l-1}\frac{1}{\epsilon(2^jR)}}{l}=\lim\limits_{j\rightarrow \infty}\frac{1}{\epsilon(2^jR)}=+\infty. \nonumber
\end{equation}
Meanwhile, we know that $0<\delta<1$ for any $R$. Therefore, letting $l\rightarrow \infty$ in \eqref{equuu3.3}, we conclude that $I_R=0$ for all sufficiently large $R$. Therefore $u$ is constant.

\qed

\section{  Eigenvalue estimates }
In this section, we derive lower bound estimations of eigenvalues of the Laplace-Beltrami operator  $\Delta$ on closed Riemannian manifolds with Bakry-\'Emery Ricci curvature bounded below and bounded the potential function.

Denote the eigenvalues of $\Delta$ by $0=\lambda_0<\lambda_1\leq \lambda_2\leq \cdots\leq \lambda_k\leq \cdots$. First, we bound $\lambda_1$ from below. According to \cite{Cheeger}, it suffices to bound Cheeger's isoperimetric constant from below.  Let us recall the definitions of isoperimetric constants. We adapt the notations and definitions in \cite{Li2}.

\begin{definition}
	Let $(M^n,g)$ be a compact Riemannian manifold (with or without boundary). For $\alpha>0$, The Neumann $\alpha$-isoperimetric constant of M is defined by
	$$	IN_\alpha(M)=\inf_{\substack{\partial\Omega_1=H=\partial\Omega_2 \\ M=\Omega_1\cup H \cup \Omega_2}}\frac{\vol(H)}{\min\{{\vol(\Omega_1),\vol(\Omega_2)}\}^{\frac{1}{\alpha}}},$$
	\\
	where the infimum is taken over all hypersurfaces $H$ dividing  $M$ into two parts, denoted by $\Omega_1$ and $\Omega_2$.
\end{definition}

In \cite{Cheeger}, Cheeger showed that
\begin{equation}\label{Cheeger 1}\lambda_1 \ge \frac{IN_1(M)^2}{4}\end{equation}
on closed manifolds. Thus, one can get a lower bound of $\lambda_1$ by bounding $IN_1(M)$ from below. We will accomplish this goal by following the method of Dai-Wei-Zhang in \cite{DWZ} and using Theorem \ref{volume element comparison}.
\begin{theorem}\label{isoperimetric estimate 1}
	Let $(M^n,g)$ be a complete Riemannian manifold with $\Ric_f\ge-Kg$ and $|f|\le L$ on $M$, where $K,\ L$ are non-negative constants. Let $\Omega$ be a bounded convex domain in $M$. Then
	for $1\le\alpha \le \frac{n}{n-1},$ we have
	\begin{equation}
		IN_\alpha(\Omega)\ge d^{-1}2^{-2n-1}5^{-n}e^{-(18-\frac{6}{\alpha})L-(17\frac{1}{6}-\frac{1}{6\alpha})Kd^2} \vol(\Omega)^{1-\frac{1}{\alpha}},
	\end{equation}
	and for $0<\alpha<1,$ we have
	\begin{equation}
		IN_{\alpha}(\Omega)\ge  d^{-1}2^{-2n-1}5^{-n}e^{-12L-17Kd^2}\vol(\Omega)^{1-\frac{1}{\alpha}},
	\end{equation}
	where $d=diam(\Omega)$, the diameter of $\Omega$.
	In particular, if $M$ is closed, then
	\begin{equation}
		IN_1(M)\ge D^{-1}2^{-2n-1}5^{-n}e^{-12L-17K D^2},
	\end{equation}
	where $D$ is an upper bound of the diameter of $M$.
	
\end{theorem}

Before starting the proof of Theorem \ref{isoperimetric estimate 1},
we need the following lemma by Gromov.

\begin{lemma}[\cite{Gro}]\label{Gromov1}
	Let $(M^n,g)$ be a complete Riemannian manifold.    Let $\Omega$ be a convex domain in $M$ and $H$  a hypersurface dividing $\Omega$	into two parts $\Omega_1,\Omega_2$. For any Borel subsets $W_i \subset \Omega_i,i=1,2$, there exists $x_1$
	in one of $W_i$, say $W_1$, and a subset $W$ in the other part $W_2$, such that
	\begin{equation}
		\vol(W) \ge \frac{1}{2}\vol(W_2),
	\end{equation}
	and for any $x_2\in W$, there is a unique minimal geodesic between  $x_1$ and $x_2$ which intersects $H$ at some $z$ with
	\begin{equation}
		d(x_1,z)\ge d(x_2,z),
	\end{equation}
	where $d(x_1,z)$ denotes the distance between $x_1$ and $z$.	
\end{lemma}
Combining  Theorem \ref{volume element comparison} and Lemma \ref{Gromov1}, we get
\begin{lemma}\label{wl2}
	Let $H,W$ and $x_1$ be as in Lemma \ref{Gromov1}. Then
	\begin{equation}
		\vol(W)\le D_12^{n-1}e^{6L+\frac{K}{2}D_1^2}\vol(H^{'}),
	\end{equation}
	where $D_1=\sup_{x\in W} d(x_1,x)$, and $H^{'}$ is the set of intersection points with $H$ of geodesics $\gamma_{x_1,x} $ for all $x \in W$.
\end{lemma}

\proof Let $\Gamma \subset S_{x_1}(M) $ be the subset of unit vectors $\theta$ such that $\gamma_\theta =\gamma_{x_1,x_2}$ for some $x_2\in W$. Set polar coordinates at $x_1$. The volume element of the metric $g$ is written as $dv=J(t,\theta,x_1)dtd\theta $ in the polar coordinates $(\theta,t) \in S_{x_1}(M) \times \mathbb{R^{+}}$. For any $\theta \in \Gamma$, let $r(\theta)$ be the radius such that $exp_{x_1}(r(\theta)\theta)\in H$. Then form Lemma \ref{Gromov1} $W\subset \{exp_{x_1}(r\theta)|\theta \in \Gamma,r(\theta) \le r \le 2r(\theta)\}$. We conclude
\[\vol(W) \le \int_{\Gamma}\int_{r(\theta)}^{2r(\theta)} J(t,\theta,x_1)dtd\theta.\]
For $r(\theta) \le t \le 2r(\theta) \le 2D_1$, by Theorem \ref{volume element comparison}, it implies
\[\frac{J(t,\theta,x_1)}{t^{n-1}}\le e^{\frac{K}{6}(t^2-r(\theta)^2)+6L} \frac{J(r(\theta),\theta,x_1)}{r(\theta)^{n-1}},\]
so \[J(t,\theta,x_1) \le e^{\frac{K}{2}D_1^2+6L}2^{n-1}J(r(\theta),\theta,x_1).\]
It gives  \[\vol(W)\le e^{\frac{K}{2}D_1^2+6L}2^{n-1}\int_{\Gamma} r(\theta)J(r(\theta),\theta,x_1)d\theta \le D_12^{n-1} e^{\frac{K}{2}D_1^2+6L}  \vol(H^{'}). \] \qed\\
Form Lemmas \ref{Gromov1} and \ref{wl2}, we immediately have

\begin{corollary}
	Let $H$ be any hypersurface dividing a convex domain $\Omega$ into two parts $\Omega_1,\Omega_2$. For any ball $B=B_r(x)$ in $M$, we have
	\begin{equation}
		\min(\vol(B\cap\Omega_1),\vol(B\cap\Omega_2)) \le  2^{n+1}re^{\frac{K}{2}d^2+6L}\vol(H\cap(B_{2r}(x))),
	\end{equation}
	where $d=diam(\Omega)$. In particular, if $B\cap\Omega$ is divided equally by $H$, then
	\begin{equation} \label{RV1}
		\vol(B\cap\Omega) \le  2^{n+2}re^{\frac{K}{2}d^2+6L}\vol(H\cap(B_{2r}(x))).
	\end{equation}
	\proof Put $W_i=B\cap\Omega_i$ in the above lemma and use $D_1 \le 2r$ and $H^{'}\subset H\cap B_{2r}(x).$ \qed\\
\end{corollary}
Now we are ready to prove Theorem \ref{isoperimetric estimate 1}.\\

\noindent{\it Proof of Theorem \ref{isoperimetric estimate 1}.} Let $H$ be any hypersurface dividing $\Omega$ into two parts, $\Omega_1$ and $\Omega_2$. We may assume that $\vol(\Omega_1) \le \vol(\Omega_2)$. For any $x\in\Omega_1$, let $r_x$ be the smallest radius such that
\[\vol(B_{r_x}(x)\cap\Omega_1)=\vol(B_{r_x}(x)\cap\Omega_2)=\frac{1}{2}\vol(B_{r_x}(x)\cap\Omega).\]
Let $d=diam(\Omega)$. By \eqref{RV1} we have,
\begin{equation}\label{CC1}
	\vol(B_{r_x}(x)\cap\Omega) \le  2^{n+2}r_xe^{\frac{K}{2}d^2+6L}\vol(H\cap(B_{2r_x}(x))).
\end{equation}
The domain $\Omega_1$ has a covering
\[\Omega_1\subset \cup_{x\in\Omega_1}B_{2r_x}(x).\]
By Vitali Covering Lemma, we can choose a countable family of disjoint balls $B_i=B_{2r_{x_i}}(x_i)$ such that $\cup_iB_{10r_{x_i}}(x_i) \supset \Omega_1.$
So\[\vol(\Omega_1)\le \sum_i \vol(B_{10r_{x_i}}(x_i)\cap\Omega_1).\]
Applying the volume comparison Theorem \ref{volume element comparison} in $\Omega_1$ gives
\[\frac{\vol(B_{10r_{x_i}}(x_i)\cap\Omega_1)}{(10r_{x_i})^n}\le e^{\frac{33}{2}Kr_{x_i}^2+6L}\frac{\vol(B_{r_{x_i}}(x_i)\cap\Omega_1)}{(r_{x_i})^n}.\]
On the other hand, since $\vol(\Omega_1) \le \vol(\Omega_2) $, we have $r_x \le d$ for any $x\in \Omega_1$. Thus,
\begin{align}
	\vol(B_{10r_{x_i}} (x_i) \cap\Omega_1)
	&\le 10^ne^{\frac{33}{2}Kd^2+6L} \vol(B_{r_{x_i}}(x_i)\cap\Omega_1)\nonumber\\ &=2^{-1}10^ne^{\frac{33}{2}Kd^2+6L} \vol(B_{r_{x_i}}(x_i)\cap\Omega)\nonumber.
\end{align}
\\
Therefore,
\begin{equation}\label{wo1}
	\vol(\Omega_1) \le 2^{-1}10^ne^{\frac{33}{2}Kd^2+6L} \sum_i\vol(B_{r_{x_i}}(x_i)\cap\Omega) .
\end{equation}
Moreover, since the balls $B_i$ are disjoint, \eqref{CC1} gives
\begin{equation}\label{wo2}
	\vol(H)\ge \sum_i \vol(B_i \cap H) \ge 2^{-n-2}e^{-\frac{K}{2}d^2-6L} \sum_i r_{x_i}^{-1} \vol(B_{r_{x_i}}(x_i)\cap\Omega).
\end{equation}
Firstly,  for $1\le \alpha \le \frac{n}{n-1}$, it follows from \eqref{wo1} and \eqref{wo2} that
\[
\begin{split}
	\frac{\vol(H)}{\vol(\Omega_1)^{\frac{1}{\alpha}}}   &\ge  \frac{2^{-n-2}e^{-\frac{K}{2}d^2-6L}}{(2^{-1}10^ne^{\frac{33}{2}Kd^2+6L})^{\frac{1}{\alpha}}}\frac{\sum_i r_{x_i}^{-1}\vol(B_{r_{x_i}}(x_i)\cap\Omega)}{(\sum_i\vol(B_{r_{x_i}}(x_i)\cap\Omega))^{\frac{1}{\alpha}}}  \\
	& \ge \frac{2^{-n-2}e^{-\frac{K}{2}d^2-6L}}{2^{-1}10^ne^{\frac{33}{2}Kd^2+6L}}\frac{\sum_ir_{x_i}^{-1} \vol(B_{r_{x_i}}(x_i)\cap\Omega)}{\sum_i\vol(B_{r_{x_i}}(x_i)\cap\Omega)^{\frac{1}{\alpha}}}  \\
	& \ge2^{-2n-1}5^{-n}e^{-12L-17Kd^2}\inf_i\frac{r_{x_i}^{-1}\vol(B_{r_{x_i}}(x_i)\cap\Omega)}{\vol(B_{r_{x_i}}(x_i)\cap\Omega)^\frac{1}{\alpha}} \\
	&= 2^{-2n-1}5^{-n}e^{-12L-17Kd^2}\inf_i r_{x_i}^{-1}\vol(B_{r_{x_i}}(x_i)\cap\Omega)^{1-\frac{1}{\alpha}}.
\end{split}
\]
We apply the volume comparison Theorem \ref{volume element comparison} in $\Omega$, then
\[\frac{\vol(B_d(x_i)\cap\Omega)}{d^n}\le e^{\frac{K}{6}(d^2-r_{x_i}^2)+6L}\frac{\vol(B_{r_{x_i}}(x_i)\cap\Omega)} {r_{x_i}^n}.\]
By $1-\frac{1}{\alpha} \ge 0$, and $n(1-\frac{1}{\alpha})-1\le0$,
we can derive \[\inf_i r_{x_i}^{-1} \vol(B_{r_{x_i}}(x_i)\cap\Omega)^{1-\frac{1}{\alpha}} \ge d^{-1}e^{-(\frac{K}{6}d^2+6L)(1-\frac{1}{\alpha})}\vol(\Omega)^{1-\frac{1}{\alpha}}. \]
Taking infimum over $H$, we conclude the following
\[IN_{\alpha}(\Omega)\ge  d^{-1}2^{-2n-1}5^{-n}e^{-(18-\frac{6}{\alpha})L-(17\frac{1}{6}-\frac{1}{6\alpha})K d^2}\vol(\Omega)^{1-\frac{1}{\alpha}}.\]
On the other hand, for $0<\alpha<1$, we have
\[
\begin{split}
	\frac{\vol(H)}{\vol(\Omega_1)^{\frac{1}{\alpha}}}&=\frac{\vol(H)}{\vol(\Omega_1) \vol(\Omega_1)^{\frac{1}{\alpha}-1}}    \ge \frac{\vol(H)}{\vol(\Omega_1) \vol(\Omega)^{\frac{1}{\alpha}-1}} \\
	& \ge \frac{2^{-n-2}e^{-\frac{K}{2}d^2-6L}}{2^{-1}10^ne^{\frac{33}{2}Kd^2+6L}}\frac{\sum_i r_{x_i}^{-1}\vol(B_{r_{x_i}}(x_i)\cap\Omega)}{\sum_i\vol(B_{r_{x_i}}(x_i)\cap\Omega)} \vol(\Omega)^{1-\frac{1}{\alpha}} \\
	& \ge d^{-1}2^{-2n-1}5^{-n}e^{-12L-17Kd^2}\vol(\Omega)^{1-\frac{1}{\alpha}} .
\end{split}
\]
\\
Taking infimum over $H$ finishes the proof. \\ \qed

From  \eqref{Cheeger 1} and Theorem \ref{isoperimetric estimate 1}, we immediately have the estimate of the first eigenvalue.
\begin{theorem}\label{thm 11}
	Let $(M^n,g)$ be a closed Riemannian manifold with $\Ric_f\ge-Kg$ and $|f|\le L$ on $M$, where $K,\ L$ are non-negative constants.
	Then
	\begin{equation}\label{the lower bound 1}
		\lambda_1\ge D^{-2}2^{-4n-4}5^{-2n}e^{-24L-34KD^2}:=\alpha_0,
	\end{equation}
	where  $D$ is an upper bound of the diameter of $M$.
\end{theorem}

Next, we derive lower bounds for $\lambda_k\ (k\ge2) $ by using the upper bound estimate of the heat kernel in Theorem \ref{Gaussian} and an argument of Li-Yau \cite{LY}.

\noindent{\it Proof of Theorem \ref{Eigenvalue Estimates theorem}.}
In Theorem \ref{Gaussian}, letting $\epsilon=1$, we have
\begin{equation}
	H(x,x,t)\le \frac{\overline{C}_1(n)e^{\overline{C}_2(n)(Kt+L)}}{\vol( B_{\sqrt{t}}(x))}\nonumber
\end{equation}		
for all $x\in M$ and $t>0$. Note that the heat kernel can be written as 	

\begin{equation}
	H(x,y,t)=\sum_{i=0}^{\infty}e^{-\lambda_it}\phi_i(x)\phi_i(y),\nonumber
\end{equation}		
where $\phi_i$ is the eigenfunction of $\Delta$ corresponding to $\lambda_i$ and $\{\phi_i \}_{i=0}^{\infty}$ form an orthonormal basis with respect to the $L^2$-norm. So we have	
\begin{equation}
	\int_{M}H(x,x,t)dx=\sum_{i=0}^{\infty}e^{-\lambda_it}\le \overline{C}_1(n)e^{\overline{C}_2(n)(Kt+L)}\int_{M}\vol( B_{\sqrt{t}}(x))^{-1}dx.\nonumber
\end{equation}		
When $\sqrt{t}\ge D$, $\vol( B_{\sqrt{t}}(x))=\vol (M)$. On the other hand, when $\sqrt{t}\le D$, since $| f|\le L$ on $M$, by the volume comparison theorem in \eqref{VC}, we get	
\begin{equation}	
	\frac{\vol( B_D(x))}{\vol( B_{\sqrt{t}}(x))}\le e^{\left(\frac{K}{6}D^2+6L\right)}\left(\frac{D}{\sqrt{t}}\right)^n\nonumber
\end{equation}	
for all $x\in M$.
Then we conclude that

\begin{equation}
	\begin{aligned}
		\sum_{i=0}^{\infty}e^{-\lambda_it}&\le \overline{C}_1(n)e^{\overline{C}_2(n)(Kt+L)}
		\begin{cases}
			e^{\left(\frac{K}{6}D^2+6L\right)}\left(\frac{D}{\sqrt{t}}\right)^n & \text{ $t\le D^2$ }, \\
			1& \text{ $t\ge D^2$} ,
		\end{cases}\\
	& \le   \begin{cases}
			\overline{C}_3(n)e^{\overline{C}_4(n)(KD^2+L)}\left(\frac{D}{\sqrt{t}}\right)^n & \text{ $t\le D^2$ }, \\
		\overline{C}_3(n)e^{\overline{C}_4(n)(Kt+L)}& \text{ $t\ge D^2$} .\nonumber
	\end{cases}:=q(t)
	\end{aligned}
\end{equation}
Fixing $k\ge1$, and taking the first $(k+1)$ terms, we get
\begin{equation}
	(k+1)e^{-\lambda_kt}\le q(t),\nonumber
\end{equation}	
that is, $k+1\le q(t)e^{\lambda_kt}$, for any $t>0$.
It is easy to see that $q(t)e^{\lambda_kt}$ is continuous for $t>0$ and
\begin{equation}
	\inf\limits_{t>0}\left(\frac{D}{\sqrt{t}}\right)^ne^{\lambda_kt}=\left(\frac{D}{\sqrt{t}}\right)^ne^{\lambda_kt}\bigg|_{t=\frac{n}{2\lambda_k}}=\left(\frac{2e}{n}\right)^{\frac{n}{2}}\cdot \left(D\sqrt{\lambda_k}\right)^n.\nonumber
\end{equation}
When $\frac{n}{2\lambda_k}\le D^2$, i.e., $\lambda_k\ge\frac{n}{2D^2}$, we have
\begin{equation}
	\inf\limits_{t>0}q(t)e^{\lambda_kt}=q(t)e^{\lambda_kt}\bigg|_{t=\frac{n}{2\lambda_k}}=	\overline{C}_3(n)e^{\overline{C}_4(n)(KD^2+L)}\left(\frac{2e}{n}\right)^{\frac{n}{2}}\cdot \left(D\sqrt{\lambda_k}\right)^n.\nonumber
\end{equation}
Hence
\begin{equation}
	\overline{C}_3(n)e^{\overline{C}_4(n)(KD^2+L)}\left(\frac{2e}{n}\right)^{\frac{n}{2}}\cdot \left(D\sqrt{\lambda_k}\right)^n\ge k+1,\nonumber
\end{equation}
that is,
\begin{equation}\label{eq6.11}
	\lambda_k\ge  \frac{\overline{C}_5(n)(k+1)^{\frac{2}{n}}}{D^2}e^{-\overline{C}_6(n)(KD^2+L)}.
\end{equation}
When $\frac{n}{2\lambda_k}\ge D^2$, i.e., $\lambda_kD^2\le\frac{n}{2}$, we have
\begin{equation}
	\inf\limits_{t>0}q(t)e^{\lambda_kt}=q(t)e^{\lambda_kt}\bigg|_{t=D^2}=	\overline{C}_3(n)e^{\overline{C}_4(n)(KD^2+L)}e^{\lambda_kD^2}.\nonumber
\end{equation}
Hence
\begin{equation}\label{c_01}
	\alpha_1:=\overline{C}_3(n)e^{\overline{C}_4(n)(KD^2+L)}e^{\frac{n}{2}}\ge k+1,
\end{equation}
which shows that there are only finitely many $\lambda_k$  in this case (the total number $k$ only depends on $n,K,L,D$). Combining \eqref{c_01} and the lower bound of $\lambda_1$ in \eqref{the lower bound 1}, for these finitely many $\lambda_k$'s, one can choose $\overline{C}_7(n)$ and $\overline{C}_{8}(n)$ such that
\begin{equation}\label{eq6.24}
	\frac{\lambda_k  D^2}{(k+1)^{\frac{2}{n}}}\ge	\frac{\lambda_1  D^2}{(k+1)^{\frac{2}{n}}}\ge 	\frac{\alpha_0  D^2}{\alpha_1^{\frac{2}{n}}}
	\ge \overline{C}_7(n)e^{-\overline{C}_{8}(n)(KD^2+L)}.
\end{equation}
Combining \eqref{eq6.11} and \eqref{eq6.24} finishes the proof.
\qed \\

Finally, we study the bottom spectrum of the Beltrami Laplacian $\Delta$ on complete Riemannian manifolds with Bakry-\'Emery Ricci curvature bounded below. Recall the definition of $\mu_1(\Delta):= \inf Spec (\Delta)$, the bottom spectrum of $\Delta$. By the variational principle, we have 
\begin{equation}\label{e7.7}
	\mu_1(\Delta)=\inf \limits_{\phi \in T\backslash \{0\} }\frac{\int_M |\nabla \phi|^2}{\int_M \phi^2},
\end{equation}
where $T$ is any class of functions such that
\begin{equation}
	C_0^\infty(M)\subset T \subset W_0^1(M).\nonumber
\end{equation}
Here $W_0^1(M)$ is the closure of $C_0^\infty(M)$ in $W^1(M)$, the $L^2$ Sobolev space.

Following Munteanu and Wang's method \cite{MW1}, the volume growth estimate and the variational principle immediately imply the upper bound of $\mu_1(\Delta)$. However, it is required that the volume should not exceed the exponential linear growth, so we cannot use Theorem \ref{volume element comparison} to obtain the growth estimate of the volume, where the calculated volume is the exponential square growth. Therefore we adopt the method of Lemma 2.1 in \cite{MW1} to obtain the volume growth estimate.
\begin{lemma}\label{volume growth estimate}
	Let $(M^n,g)$ be a complete Riemannian manifold with $\Ric_f\ge -(n-1)Kg$ for some constant $K\ge0$. Assume that there exist non-negative constants $\tilde{a}$ and $\tilde{b}$ such that 
	\begin{equation}
		|f|(x)\le \tilde{a}r(x,o)+\tilde{b}\ for\ all\ x\in M.\nonumber
	\end{equation}
Then for any $\epsilon>0$, there exists a constant $\tilde{A}(\epsilon)>0$ such that the volume upper bound 
\begin{equation}
	\vol (B_R(o))\le \tilde{A}(\epsilon)e^{(2\tilde{a}+(n-1)(\sqrt{K}+\epsilon))R}\nonumber
\end{equation}
holds for all $R>0$. Here, $o$ is a fixed point on $M$, and $r(x,o)$ denotes the distance from $x$ to $o$.
\proof By the Bochner formula, we have for $r(x)=r(x,o)$
\begin{equation}
	\begin{aligned}
		0=\frac{1}{2}\Delta|\nabla r|^2&=|\Hess r|^2+\left<\nabla\Delta r,\nabla r\right>+\Ric(\partial r,\partial r) \\
		&\ge \frac{(\Delta r)^2}{n-1}+\partial_r(\Delta r)+\Ric (\partial_r,\partial_r)\\
		&\ge \frac{(\Delta r)^2}{n-1}+\partial_r(\Delta r)-f''(r)-(n-1)K.\nonumber
		\end{aligned}
\end{equation}
Integrating this inequality from $1$ to $r$, we get 
\begin{equation}
	\frac{1}{n-1}\int_{1}^{r}(\Delta t)^2dt +\Delta r-f'(r)\le (n-1)Kr+\tilde{b}_0\nonumber
\end{equation}
for some constant $\tilde{b}_0>0$  independent of $r$.\\
Then for any $r\ge1$,
\begin{equation}\label{s1}
	\Delta_f r+\frac{1}{n-1}\int_{1}^{r}\left(\Delta_f t+f'(t)\right)^2dt\le (n-1)Kr+\tilde{b}_0,
\end{equation}
The Cauchy-Schwarz inequality implies that 
\begin{equation}
	\int_{1}^{r}\left(\Delta_f t+f'(t)\right)^2dt\ge\frac{1}{r-1}\left(\int_{1}^{r}\left(\Delta_f t+f'(t)\right)dt\right)^2.\nonumber
\end{equation}
Therefore, from \eqref{s1} we obtain
\begin{equation}\label{s2}
	\Delta_f r+\frac{1}{(n-1)r}\left(f(r)-f(1)+\int_{1}^{r}\left(\Delta_f t\right)dt\right)^2\le (n-1)Kr+\tilde{b}_0.
\end{equation}
We now claim that for any $r\ge1$ and any $\epsilon>0$,
\begin{equation}\label{s3}
	\int_{1}^{r}\left(\Delta_f t\right)dt\le (\tilde{a}+(n-1)(\sqrt{K}+\epsilon))r+\tilde{a}+2\tilde{b}+\frac{\tilde{b}_0}{\sqrt{K}+\epsilon}.
\end{equation}
To prove this, define
\begin{equation}
	v(r):=(\tilde{a}+(n-1)(\sqrt{K}+\epsilon))r+\tilde{a}+2\tilde{b}+\frac{\tilde{b}_0}{\sqrt{K}+\epsilon}-\int_{1}^{r}\left(\Delta_f t\right)dt.\nonumber
\end{equation}
We show instead that $v(r)>0$ for all $r\ge1$. Clearly, $v(1)>0$.\\
Suppose that $v$ does not remain positive for all $r\ge1$ and let $R>1$ be the first number such that $v(R)=0$. Then $v'(R)\le 0$ and 
\begin{equation}
	\int_{1}^{R}\left(\Delta_f t\right)dt=(\tilde{a}+(n-1)(\sqrt{K}+\epsilon))R+\tilde{a}+2\tilde{b}+\frac{\tilde{b}_0}{\sqrt{K}+\epsilon}.\nonumber
\end{equation}
In other words,
\begin{equation}
	\begin{aligned}
		&\frac{1}{(n-1)R}\left(f(R)-f(1)+\int_{1}^{R}\left(\Delta_f t\right)dt\right)^2\\
		&=\frac{1}{(n-1)R}\left(f(R)-f(1)+(\tilde{a}+(n-1)(\sqrt{K}+\epsilon))R+\tilde{a}+2\tilde{b}+\frac{\tilde{b}_0}{\sqrt{K}+\epsilon}\right)^2\\
		&\ge\frac{1}{(n-1)R}\left((n-1)(\sqrt{K}+\epsilon)R+\frac{\tilde{b}_0}{\sqrt{K}+\epsilon}\right)^2\\
		&\ge (n-1)(\sqrt{K}+\epsilon)^2R+2\tilde{b}_0.\nonumber
	\end{aligned}
\end{equation}
\end{lemma}
Plugging this into \eqref{s2}, we conclude that $\Delta_f R\le-\tilde{b}_0<0$, so $v'(R)=\tilde{a}+(n-1)(\sqrt{K}+\epsilon)-(\Delta_f R)>0$, which is a contradiction.\\
We have thus proved \eqref{s3} is true for any $r\ge1$ and $\epsilon>0$, so 
\begin{equation}
	\ln J(r)-\ln J(1)\le (2\tilde{a}+(n-1)(\sqrt{K}+\epsilon))r+2\tilde{a}+4\tilde{b}+\frac{\tilde{b}_0}{\sqrt{K}+\epsilon}.\nonumber
\end{equation}
In particular, for $R\ge1$, any $\epsilon>0$ we have the volume bound of the form
\begin{equation}
	\vol (B_R(o)) \le \tilde{b}_1e^{(2\tilde{a}+(n-1)(\sqrt{K}+\epsilon))R},\nonumber
\end{equation}
where the constant $\tilde{b}_1$ depends on $\tilde{a},\tilde{b},\vol (B_1(o)),\epsilon,K$ and $n$. \\ \qed

\noindent{\it Proof of Theorem \ref{upper bound of spectrum}.} 
Let $R>1$ and $\psi$ a cut-off function on $B_R(o)$ such that $\psi=1$ on $B_{R-1}(o)$ and $|\nabla \psi|\le 2$. Set $\phi(y):=e^{-\frac{(2\tilde{a}+(n-1)(\sqrt{K}+\epsilon)+\delta)}{2}r(y,o)}\psi(y)$ as a test function in the variational principle \eqref{e7.7} for $\mu_1(\Delta)$, where $\delta>0$ and $\epsilon>0$ are arbitrary positive constants. Then by Lemma \ref{volume growth estimate}, we obtain
\begin{equation}
	\mu_1(\Delta) \le \frac{1}{4}\left(2\tilde{a}+(n-1)(\sqrt{K}+\epsilon)+\delta\right)^2.\nonumber
\end{equation}
Since $\epsilon$ and $\delta$ are arbitrary, then $\mu_1(\Delta) \le \frac{1}{4}\left(2\tilde{a}+(n-1)\sqrt{K}\right)^2.$ 

In the case that $f$ is of sublinear growth, we can take $\tilde{a}=0$. Therefore, $\mu_1(\Delta)\le\frac{1}{4}(n-1)^2K$ and the theorem is proved. \\ \qed

\section*{Data Availability}
No data was used for the research described in the article.

\section*{Acknowledgements}

Research is partially supported by NSFC Grant No. 11971168, Shanghai Science and Technology Innovation Program Basic Research Project STCSM 20JC1412900, and Science and Technology Commission of Shanghai Municipality (STCSM) No. 22DZ2229014.


\begin{thebibliography}{99}











\bibitem{BE}  D. Bakry and M. \'Emery, \emph{Diffusions hypercontractives.}  (French) [Hypercontractive diffusions] S\'eminaire de probabilit\'es, XIX, 1983/84, 177-206, Lecture Notes in Math., 1123, Springer, Berlin, 1985.
\bibitem{BQ}  D. Bakry and Z. Qian, {\it Volume comparison theorems without Jacobi fields}, Current trends in potential theory, 115-122, Theta Ser. Adv. Math., 4, Theta, Bucharest, 2005.
\bibitem{Bu} P. Buser, \emph{A note on the isoperimetric constant}, Ann. Sci. Ecole Norm. Sup. 15 (1982), 213-230
\bibitem{CSY} S.Y. Cheng, \emph{Eigenvalue comparison theorems and its geometric applications}, Math. Z. 143 (1975), no. 3, 289–297.
\bibitem{CZ} H.-D. Cao and D. Zhou, {\it On complete gradient shrinking Ricci solitons}, J. Diffeential Geom., 85 (2010), 175-186.
\bibitem{CG} G. Carron, {\it Harmonic functions on manifolds whose large spheres are small}, Ann. Math. Blaise Pascal 23 (2016), no. 2, 249-261.
\bibitem{Cheeger}   J. Cheeger, {\it A lower bound for the smallest eigenvalue of the Laplacian}, Problems in analysis (Papers dedicated to Salomon Bochner, 1969), pp. 195-199. Princeton Univ. Press, Princeton, N. J., 1970.
\bibitem{Chen} B.-L. Chen, {\it Strong uniqueness of the Ricci flow}, J. Differential Geom. 82(2): 363-382.
\bibitem{DWZ} X. Dai, G. Wei and Z. Zhang, {\it Neumann isoperimetric constant estimate for convex domains}, Proc. Amer. Math. Soc. 146, no. 8, 3509-3514.
\bibitem{Da} E.B. Davies, \emph{Heat Kernels and Spectral Theory}, Cambridge Tracts in Mathematics, vol. 92, Cambridge University Press, Cambridge ,1989
\bibitem{FLZ}  F. Fang, X.-D. Li and Z. Zhang, {\it Two generalizations of Cheeger-Gromoll splitting theorem via Bakry-Emery Ricci curvature}, Ann. Inst. Fourier (Grenoble) 59 (2009), no. 2, 563-573.
\bibitem{GKW} G. J. Galloway, M. A. Khuri and E. Woolgar, {\it A Bakry-\'Emery almost splitting result with applications to the topology of black holes}, Comm. Math. Phys. 384 (2021), no. 3, 2067–2101.
\bibitem{Gr} A. Grigor'yan, \emph{Heat kernel and analysis on manifolds}, AMS/IP Studies in Advanced Mathematics, 47. American Mathematical Society, Providence, RI; International Press, Boston, MA, 2009.
\bibitem{Gr2} A. Grigor'yan, {\it Upper bounds of derivatives of the heat kernel on an arbitrary complete manifold}, J. Funct. Anal. , 127 no.2, (1995) 363-389.
\bibitem{Gro} M. Gromov, {\it Paul Levy's isoperimetric inequality}, Preprint IHES, 1979.
\bibitem{GPSS} B. Guo, D. H. Phong, J. Song and J. Sturm, {\it Compactness of K\"ahler-Ricci solitons on Fano manifolds,} Pure Appl. Math. Q. 18 (2022), no. 1, 305-316.
\bibitem{Ham} R. Hamilton, {\it The formation of singularities in the Ricci flow}, Surveys in differential geometry, Vol. II (Cambridge, MA, 1993), 7-136, Int. Press, Cambridge, MA, 1995.
\bibitem{LLW} H. Li, Y. Li and B. Wang, {\it On the structure of Ricci shrinkers}, J. Funct. Anal. 280 (2021), no. 9, Paper No. 108955, 75 pp.
\bibitem{Li} P. Li, \emph{Uniqueness of $L^1$ solutions for the Laplace equation and the heat equation on Riemannian manifolds}, J. Differ. Geom. 20 (1984), 447-457.
\bibitem{Li2} P. Li, {\it Geometric Analysis}, Cambridge Studies in Advanced Mathematics, vol. 134, Cambridge University Press, Cambridge, 2012.
\bibitem{LiSc} P. Li and R. Schoen, {\it $L^p$ and mean value properties of subharmonic functions on Riemannian manifolds}, Acta Math. 153 (1984), no. 3-4, 279-301.
\bibitem{LY} P. Li and S.-T. Yau, \emph{On the parabolic kernel of the Schrodinger operator}, Acta Math. 156 (1986), 153-201.
\bibitem{Lix1}  X.-D. Li, {\it Liouville theorems for symmetric diffusion operators on complete Riemannian manifolds}, J. Math. Pures Appl. 84, 1295-1361 (2005)
\bibitem{Lix2} X.-D. Li, {\it Perelman's entropy formula for the Witten Laplacian on Riemannian manifolds via Bakry-Emery Ricci curvature}, Math. Ann. 353 (2012), no. 2, 403-437.
\bibitem{Liy}  Y. Li, {\it Li-Yau-Hamilton estimates and Bakry-Emery-Ricci curvature}, Nonlinear Anal. 113 (2015), 1-32.
\bibitem{Liu} G. Liu, {\it Stable weighted minimal surfaces in manifolds with non-negative Bakry-Emery Ricci tensor}, Comm. Anal. Geom. 21 (2013), no. 5, 1061-1079.
\bibitem{Lo} J. Lott, {\it Some geometric properties of the Bakry-\'Emery-Ricci tensor}, Comment. Math. Helv. 78 (2003), no. 4, 865-883.
\bibitem{LV}  J. Lott and C. Villani, \emph{Ricci curvature for metric-measure spaces via optimal transport.} Ann. of Math. (2) 169 (2009), no. 3, 903-991.
\bibitem{MuSe} O. Munteanu and N. Sesum, {\it On gradient Ricci solitons}, J. Geom. Anal. 23 (2013), no. 2, 539-561.
\bibitem{MW} O. Munteanu and J. Wang, \emph{Smooth metric measure spaces with nonnegative curvature}, Commun. Anal.Geom. 19 (2011), 451-486.
\bibitem{MW1} O. Munteanu and J. Wang, \emph{Analysis of weighted Laplacian and applications to Ricci solitons}. Comm. Anal. Geom. 20 (2012), no. 1, 55–94.
\bibitem{Per} G. Perelman, \emph{The entropy formula for the Ricci flow and its geometric applications,} arXiv:math/0211159.
\bibitem{PRS} S. Pigola, M. Rimoldi, A. Setti, {\it Remarks on non-compact gradient Ricci solitons}, Math. Z. 268 (2011), 777-790.
\bibitem{Qian} Z. Qian, {\it A comparison theorem for an elliptic operator}, Potential Anal. 8 (1998), no. 2, 137-142.
\bibitem{SC} L. Saloff-Coste, \emph{Aspects of Sobolev-Type Inequalities}, London Math. Soc. Lecture Note Ser., vol. 289, Cambridge University Press, Cambridge, 2002.
\bibitem{St1}  K.-T. Sturm, \emph{On the geometry of metric measure spaces.I.} Acta Math. 196 (2006), no. 1, 65-131.
\bibitem{St2}  K.-T. Sturm, \emph{On the geometry of metric measure spaces. II.} Acta Math. 196 (2006), no. 1, 133-177.
\bibitem{SuZh} Y.-H. Su and H.-C. Zhang, {\it Rigidity of manifolds with Bakry-\'Emery Ricci curvature bounded below}, Geom. Dedicata 160 (2012), 321-331.
\bibitem{WZ} F. Wang and X. Zhu, {\it The structure of spaces with Bakry-\'Emery Ricci curvature bounded below}, J. Reine Angew. Math. 757 (2019), 1-50
\bibitem{WeWy} G. Wei and W. Wylie, {\it Comparison geometry for the Bakry-Emery Ricci tensor}, J. Diff. Geom. 83 (2009), no. 2, 377-405.
\bibitem{Wu1}  J.-Y. Wu, {\it Upper bounds on the first eigenvalue for a diffusion operator via Bakry-\'Emery Ricci curvature}, J. Math. Anal. Appl. 361 (2010), no. 1, 10-18.
\bibitem{Wu2} J.-Y. Wu, {\it Upper bounds on the first eigenvalue for a diffusion operator via Bakry-\'Emery Ricci curvature II}, Results Math. 63 (2013), no. 3-4, 1079-1094.
\bibitem{Wu3} J.-Y. Wu, {\it $L^p$-Liouville theorems on complete smooth metric measure spaces},
Bull. Sci. Math. 138 (2014), no. 4, 510-539.
\bibitem{WuWu} J.-Y. Wu and P. Wu, \emph{Heat kernel on smooth metric measure spaces with nonnegative curvature}, Math. Ann. 362 (2015), no. 3-4, 717-742.
\bibitem{WuWu2} J.-Y. Wu and P. Wu,  {\it Heat kernel on smooth metric measure spaces and applications,} Math. Ann. 365 (2016), no. 1-2, 309-344.
\bibitem{Yang} N. Yang, {\it A note on nonnegative Bakry-\'Emery Ricci curvature}, Arch. Math. 93 (2009), 491-496.
\bibitem{Yau1} S.-T. Yau, {\it Harmonic functions on complete Riemannian manifolds}, Comm. Pure Appl. Math. 28 (1975), 201-228.
\bibitem{Yau2} S.-T. Yau, {\it Some function-theoretic properties of complete Riemannian manifolds and their applications to geometry,} Indiana Univ. Math. J. 25 (1976), no. 7, 659-670.
\bibitem{ZZ1} Q. S. Zhang and M. Zhu, \emph{New volume comparison results and applications to degeneration of Riemannian metrics}, Adv. Math. 352 (2019), 1096-1154.
\bibitem{ZZ2} Q. S. Zhang and M. Zhu, {\it Bounds on harmonic radius and limits of manifolds with bounded Bakry-\'Emery Ricci curvature,} J. Geom. Anal. 29 (2019), no. 3, 2082-2123.


\end{thebibliography}
	\end{document}